\documentclass[a4paper, 11pt]{article}
\usepackage[latin1]{inputenc}
\usepackage[francais]{babel}
\usepackage{epsfig}
\usepackage{amsmath}
\usepackage{amssymb}
\usepackage{latexsym}

\newtheorem{thm}{Th{\'e}or{\`e}me}[section]
\newtheorem{prop}[thm]{Proposition}
\newtheorem{cor}[thm]{Corollaire}
\newtheorem{lem}[thm]{Lemme}
\newtheorem{defn}{D{\'e}finition}[section]
\newtheorem{conj}[thm]{Conjecture}

\newcommand{\ie}{{\it ie.}}
\newcommand{\cf}{{\it cf.} }
\newcommand{\me}{\hat{\mu}^{\mathrm{ess}}_{\mathcal{L}}}
\newcommand{\de}{\mathrm{deg}_{\mathcal{L}}}

\date{}
\begin{document}
\title{Une minoration du minimum essentiel sur les variétés abéliennes}

\newenvironment{dem}{\textbf{Preuve}\par}
{\begin{flushright}$\Box$\end{flushright}}

\author{\\ Aur\'elien Galateau\\}
\maketitle

\bigskip

{\small {\sc Résumé.} On étend à la codimension générale la minoration du minimum essentiel sur les variétés abéliennes obtenue dans \cite{Galateau08}, sous une conjecture concernant leurs idéaux premiers ordinaires. Cette minoration est optimale ``à $\epsilon$ près'' en le degré de la sous-variété. Une nouvelle approche de la phase de transcendance, suivant \cite{Amoroso07}, simplifie le lemme de zéros et sa combinatoire, ce qui permet de mettre en place l'argument de descente. Celui-ci est rendu beaucoup plus délicat dans le cadre abélien par l'absence, en général, d'un relèvement du morphisme de Frobenius en caractéristique nulle. 

\bigskip

{\sc Abstract.} We extend to the general codimension a lower bound for the essential minimum on abelian varieties found in \cite{Galateau08}, under a conjecture about ordinary primes in abelian varieties. This lower bound is the best expected, ``up to an $\epsilon$'', in the degree of the subvariety. Following \cite{Amoroso07}, we change the transcendance phasis in order to simplify the zero lemma and its combinatorics. The last argument is a descent procedure on varieties, which is far more intricate in the abelian setting since there is not, in general, a lifting of the Frobenius in characteristic zero.}  

\bigskip

\section{Introduction}

On poursuit dans ce travail la minoration de la hauteur normalisée sur les variétés abéliennes commencée dans \cite{Galateau08}, dont on généralise le résultat principal en codimension quelconque. 

\bigskip

Commençons par rappeler les origines de notre problème. Soit $C$ une courbe algébrique de genre $g \geq 2$ définie sur $\overline{\mathbb{Q}}$ et plongée dans sa jacobienne $J(C)$. On note $\hat{h}$ la hauteur canonique sur $J(C)$. En 1981, Bogomolov conjecture:
\begin{thm}
\label{1}
Il existe $\epsilon > 0$ tel que $\{ x \in C(\bar{\mathbb{Q}}), \hat{h}(x) \leq \epsilon \}$ est fini.
\end{thm}

Plus généralement, soit $V$ une sous-variété d'une variété abélienne munie d'un fibré $\mathcal{L}$ ample et symétrique. On note $\hat{h}_{\mathcal{L}}$ la hauteur de Néron-Tate associée à ce fibré. On dit que $V$ est {\it de torsion} si $V$ est la translatée d'une sous-variété abélienne par un point de torsion. Une courbe algébrique de torsion est en particulier de genre $1$. On définit aussi, pour $V$ un fermé de Zariski inclus dans $A$:
\begin{defn}
Le minimum essentiel de $V$ est:
\begin{displaymath}
\me(V)= \mathrm{inf} \{ \theta > 0,
\overline{V(\theta)}^Z= V(\overline{\mathbb{Q}}) \},
\end{displaymath}
\begin{center}
où $V(\theta)= \{ x \in V(\overline{\mathbb{Q}}), \hat{h}_{\mathcal{L}}(x) \leq
\theta \}$ et $\overline{V(\theta)}^Z$ est son adhérence de Zariski.
\end{center}
\end{defn}
On peut alors étendre aux sous-variétés des variétés abéliennes le théorème \ref{1}: 
\begin{thm}
\label{2}
Soit $V$ une sous-variété propre d'une variété abélienne $A$ (\ie: $V \subsetneq A$).
Le minimum essentiel de $V$ est nul si et seulement si $V$ est de torsion.
\end{thm}
Le théorème \ref{1} a été démontré par Ullmo (\cf \cite{Ullmo96}) et le théorème \ref{2} par Zhang (\cf \cite{Zhang98}), en utilisant les propriétés d'équirépartition démontrées dans \cite{Szpiro-Ullmo-Zhang97}. Le résultat analogue est vrai si on remplace $A$ par un tore (\cf \cite{Zhang92}) ou plus généralement par une variété semi-abélienne (\cf \cite{David-Philippon00}).

\bigskip

\subsection{Versions quantitatives du problème de Bogomolov}

Ces conjectures étant démontrées, on peut s'intéresser à une version quantitative, en précisant $\epsilon$ dans le théorème \ref{1} ou en 
minorant le minimum essentiel d'une variété qui n'est pas
de torsion. Grâce au théorème des minima successifs démontré par
Zhang (dans \cite{Zhang-Positive95}), ceci revient à minorer la hauteur d'une telle variété.
Depuis les travaux de Bombieri et Zannier (voir  \cite{Bombieri-Zannier95} pour le cas torique et \cite{Bombieri-Zannier96} pour le cas abélien, le théorème $1$ de cet article étant vrai pour toute variété abélienne après la démonstration par Zhang de la conjecture de Bogomolov généralisée), on sait qu'on peut espérer obtenir une borne ``uniforme'' pour le minimum essentiel, ne dépendant que du degré de $V$ et de la variété abélienne $A$.

Remarquons que pour obtenir de telles minorations, on doit exclure les translatés de sous-variétés abéliennes par des points qui ne sont pas de torsion. En effet, si $V=x+B \subsetneq A$, le minimum essentiel est relié à la hauteur de la projection de $x$ dans le quotient $A/B$ (pour plus de détails, voir \cite{Litcanu99}) et on peut le faire tendre vers $0$ en fixant le degré et la dimension de $V$. 

\bigskip

Amoroso et David obtiennent une majoration quasi-optimale en le degré
de $V$ pour les sous-variétés des tores (\cite{Amoroso-David03}). Le degré y est en fait remplacé
par un invariant plus fin et plus approprié aux techniques diophantiennes, l'indice
d'obstruction.
\begin{defn}Soit $V$ un fermé de Zariski propre de $S$ une variété semi abélienne munie d'un fibré ample.
On appelle indice d'obstruction de $V$, noté $\omega( V)$:
\begin{displaymath}
\omega( V)= \mathrm{inf} \{
\mathrm{deg}(Z) \},
\end{displaymath}
où l'infimum est pris sur l'ensemble des hypersurfaces (non nécessairement irréductibles) de $S$ contenant $V$.
\end{defn}
Dans un schéma d'approximation diophantienne, ils utilisent les propriétés de ramification des corps cyclotomiques, qu'ils traduisent en un résultat métrique leur permettant d'extrapoler une fonction auxiliaire nulle sur $V$. Par le plongement standard $\mathbb{G}_m^n \hookrightarrow \mathbb{P}^n$, on dispose d'une hauteur projective $h$ sur les points de $\mathbb{G}_m^n$, et un minimum essentiel $\hat{\mu}^{\mathrm{ess}}$ sur les sous-variétés de $\mathbb{G}_m^n$. Ils démontrent alors:
\begin{thm}
\label{3}
Soit $V$ une sous-variété propre (et irréductible) de $\mathbb{G}_{m}^{n}$ de codimension $k$ qui n'est contenue dans aucun translaté d'un sous-tore propre de $\mathbb{G}_{m}^{n}$. On a:
\begin{displaymath}
\hat{\mu}^{\mathrm{ess}}(V) \geq \frac{c(n)}{\omega(V)} \times \Big(\mathrm{log}\big(3 \omega(V)\big)\Big)^{-\lambda(k)},
\end{displaymath}
où $c(n)$ est un réel strictement positif et $\lambda(k)=(9(3k)^{(k+1)})^{k}.$
\end{thm}
Ce théorème est optimal aux termes logarithmiques près (voir la conjecture 1.2 de \cite{Amoroso-David03}, et la remarque qui suit), comme le résultat de Dobrowolski sur le problème de Lehmer et ses généralisations. 

\bigskip

Dans le cas des variétés abéliennes, on dispose déjà de résultats quantitatifs, mais la
dépendance en le degré est moins bonne. On cite ici une version faible du théorème principal de \cite{David-Philippon02}, obtenue en comparant le rayon d'injectivité à la hauteur de $A$ suivant le lemme 6.8 de cette même référence:
\begin{thm}
\label{4}
Soit $A$ une variété abélienne de genre $g \geq 2$ définie sur $\overline{\mathbb{Q}}$, principalement polarisée par un fibré $\mathcal{L}$, et $V$ une sous-variété algébrique de $A$ qui n'est pas translatée d'une sous-variété abélienne. Il existe alors une fonction explicite et strictement positive $C(g)$ telle que:
\begin{displaymath}
\me(V) \geq \frac{ C(g) }{ d^{g+1} \mathrm{max} \{1, h(A) \} ^{g+1} \de(V)^{2(g+1)^2}},
\end{displaymath}
où $d$ est le degré d'un corps de définition de $A$ et $h(A)$ est la hauteur projective de l'origine de $A$ dans le plongement associé à $\mathcal{L}^{\otimes 16}$.
\end{thm}
Cette minoration est monomiale inverse en le degré, alors que dans le cas torique, elle est linéaire inverse en l'indice d'obstruction (aux termes logarithmiques près), ce qui, en vertu d'un théorème de Chardin (voir \cite{Chardin89}, 2, Corollaire 2), correspond à une minoration en deg$(V)^{-1/\mathrm{codim}(V)}$.

Remarquons enfin que l'hypothèse du théorème \ref{4} ($V$ n'est pas un translaté de sous-variété abélienne propre) est plus faible que son analogue torique dans le théorème \ref{3} ($V$ n'est pas {\it incluse} dans un translaté de sous-tore propre); cette différence se ressent dès qu'on obtient des résultats comparables au théorème \ref{3}, et on peut préciser la minoration sous l'hypothèse faible en faisant intervenir la dimension du plus petit translaté de sous-tore propre contenant $V$ (voir le corollaire 1.6 de \cite{Amoroso-David03}). 

\subsection{Résultats et applications}

Soit $A$ une variété abélienne définie sur $K$ un corps de nombres, munie d'un fibré $\mathcal{L}$ ample et symétrique. Pour $\mathcal{A}$ un modèle entier de $A$ sur $\mathcal{O}_K$, on commence par poser l'hypothèse suivante: \\
\\ 
\textbf{Hypothèse H} Il existe une densité positive de premiers $\mathfrak{p}$ de $\mathcal{O}_K$ en lesquels le $p$-rang de la fibre spéciale $\mathcal{A}_{\mathfrak{p}}$ est égal à $g$. \\
\\ 
Sous cette hypothèse, on a démontré dans \cite{Galateau08} l'analogue du théorème \ref{3} en petite codimension ($k \leq 2$). 

\bigskip

Le but de ce travail est de généraliser à la codimension quelconque ce résultat. Suivant une idée nouvelle d'Amoroso (voir \cite{Amoroso07}), on modifie les arguments de transcendance, ce qui permet d'insérer un principe de tiroirs dans la combinatoire du lemme de zéros. Grâce à ce raffinement, on peut mettre en place la technique de descente, ce qui n'était possible qu'en codimension $2$ dans notre travail précédent. On démontre ici: 
\begin{thm}
\label{5}
On suppose qu'il existe un modèle entier $\mathcal{A}$ de $A$ vérifiant $\textbf{H}$. Soit $V$ une sous-variété propre (et irréductible)
de codimension $k$ quelconque dans $A$. Si $V$ n'est contenue dans aucun translaté
d'une sous-variété abélienne propre de $A$, on a:
\begin{displaymath}
\me(V) \geq \frac{C(A)}{\omega(V)} \times
\Big(\mathrm{log}\big(3 \de(V)\big)\Big)^{-\lambda(k)},
\end{displaymath}
où $C(A)$ est un réel strictement positif ne dépendant que de $A$ et
$\lambda(k)=(5g(k+1))^{k+1}$.
\end{thm}
\textbf{Remarque} La constante $\lambda(k)$ obtenue avec ce nouveau schéma de preuve est sensiblement meilleure que celle du théorème \ref{3}.

\bigskip

Par l'inégalité des minima successifs (voir {\it infra}, \ref{2.1} pour cette inégalité et pour la définition de $\hat{h}$ sur les sous-variétés de $A$), on en déduit:
\begin{cor}
\label{6}
S'il existe un modèle entier de $A$ vérifiant $\textbf{H}$, pour toute sous-variété $V$ propre et irréductible de codimension $k \geq 2$ dans $A$ qui n'est contenue dans aucun translaté
d'une sous-variété abélienne propre de $A$, on a:
\begin{displaymath}
\hat{h}_{\mathcal{L}}(V) \geq C'(A),
\end{displaymath}
où $C'(A)$ est un réel strictement positif ne dépendant que de $A$.
\end{cor}

La seule restriction qui empêche d'obtenir une minoration inconditionnelle est donc l'hypothèse $\textbf{H}$. Elle est en fait l'objet de la conjecture suivante (voir \cite{Pink98}, 7) :
\begin{conj} 
\label{conjecture ordinaire}
Pour toute variété abélienne $A$ définie sur $\overline{\mathbb{Q}}$, il existe un corps de définition $K$ de $A$ et un modèle entier $\mathcal{A}$ de $A$ sur $\mathcal{O}_K$ vérifiant \textbf{H}.
\end{conj}
Sous cette conjecture, la minoration du théorème \ref{5} est vérifiée pour toutes les variétés abéliennes définies sur $\overline{\mathbb{Q}}$. 
En dimension $1$, pour une courbe elliptique $E$, le résultat est connu (voir \cite{Serre68}, IV, 13 et la théorie de la multiplication complexe). La validité de \textbf{H} a été étendue aux surfaces abéliennes par les travaux de Katz et Ogus
(\cf \cite{Ogus83} 2.7, en remarquant par le théorème de Chebotarev que les premiers de degré 1
ont une densité positive). L'hypothèse \textbf{H} est encore vraie pour les variétés abéliennes à multiplication complexe. Des conditions suffisantes, portant sur les groupes de monodromie
$G_l$ (associés à chaque nombre premier $l$) de la variété abélienne, ont été données par Noot
(voir \cite{Noot95}, 2) puis Pink (\cite{Pink98}, 7).

\bigskip

Une telle minoration du minimum essentiel s'applique au problème de l'intersection d'une variété avec des sous-groupes algébriques, en direction de conjectures formulées de façon indépendante par Zilber sur les variétés semi-abéliennes, et Pink sur les variétés de Shimura mixtes (voir \cite{Zilber02} et \cite{Pink05}, conjectures 1.2 et 1.3). Pour $S$ un sous-ensemble de la variété abélienne $A$, on note:
\begin{displaymath}
S_{\epsilon}=\{ x+y, \ x \in S, \ y \in A, \ \hat{h}_{\mathcal{L}}(y) \leq \epsilon \},
\end{displaymath}
ainsi que: 
\begin{displaymath}
\mathcal{H}= \bigcup_{\mathrm{codim} H \geq 2} H,
\end{displaymath}
la réunion portant sur tous les sous-groupes algébriques de $A$ ayant la codimension prescrite. On garde les mêmes définitions avec le tore $\mathbb{G}_m^n(\overline{\mathbb{Q}})$ à la place de $A$.
En utilisant le théorème \ref{3}, Habegger a démontré le résultat suivant (voir \cite{Habegger06}):  
\begin{thm}
\label{7}
Soit $C$ une courbe dans $\mathbb{G}_m^n(\overline{\mathbb{Q}})$ qui n'est pas incluse dans le translaté d'un sous-tore propre. Alors il existe $\epsilon >0$ tel que $C \cap \mathcal{H}_{\epsilon}$ est fini.
\end{thm}
Ce théorème a permis à Maurin de démontrer la conjecture de Zilber pour une courbe plongée dans un tore (\cf \cite{Maurin07}):
\begin{thm}
\label{8}
Soit $C$ une courbe algébrique irréductible de $\mathbb{G}_m^n(\overline{\mathbb{Q}})$ non incluse dans un sous-groupe algébrique propre (éventuellement réductible). Alors $C \cap \mathcal{H}$ est fini.
\end{thm}
Il améliore ainsi le résultat principal de l'article fondateur de Bombieri, Masser et Zannier (\cf \cite{Bombieri-Masser-Zannier99}). Dans le cadre abélien, des estimations de type Lehmer sur les points de petite hauteur permettent essentiellement de traiter le cas des variétés abéliennes à multiplication complexe (voir \cite{Viada02}, \cite{Remond-Viada03}, \cite{Ratazzi05}). 

\bigskip

Soit $E$ une courbe elliptique définie sur $\overline{\mathbb{Q}}$, $g \geq 2$ un entier et $\Gamma$ un sous-groupe de $E^g$ de rang fini. Le théorème \ref{5} rend désormais inconditionnel le résultat de Viada (\cf \cite{Viada06}): 
\begin{thm}
Soit $C$ une courbe algébrique incluse dans $E^g$.
\begin{itemize}
\item Si $C$ n'est pas incluse dans un sous-groupe algébrique propre, il existe $\epsilon > 0$ tel que $C \cap \mathcal{H}_{\epsilon}$ soit fini.
\item Si $C$ n'est pas incluse dans un translaté de sous-variété abélienne propre, il existe $\epsilon > 0$ tel que $C \cap (\mathcal{H} + \Gamma_{\epsilon})$ soit fini. 
\end{itemize}
\end{thm}
Le premier point de ce théorème implique la conjecture de Zilber pour les courbes plongées dans une puissance de courbe elliptique. En raffinant son argument, Viada a depuis obtenu des résultats pour une sous-variété générale d'une puissance de courbe elliptique (voir \cite{Viada08}).

\subsection{Plan de l'article}

La preuve du théorème \ref{5} est une adaptation du nouveau schéma de transcendance utilisé par Amoroso dans \cite{Amoroso07}. Elle est écrite dans le formalisme classique de la géométrie diophantienne. La deuxième partie rassemble quelques rappels sur les variétés abéliennes. On y introduit notamment le plongement étiré, qui permet de négliger la constante de comparaison entre la hauteur projective et la hauteur de Néron-Tate. On construit aussi un bon ensemble de premiers de densité positive.

\bigskip

Dans la troisième partie, on choisit une base de dérivations adaptée, puis on démontre un lemme clé $p$-adique sur les points de torsion des variétés abéliennes. On transporte ensuite cette propriété à un polynôme en les coordonnées affines de la variétés abélienne. Dans la partie suivante, on en déduit une inégalité entre la fonction de Hilbert d'un fermé avec multiplicité, et la fonction de Hilbert de ce même fermé translaté par un sous-groupe de torsion bien choisi, sans la multiplicité. Un bon choix de paramètres nous assure alors de l'existence d'une section non-nulle, mais nulle sur ce translaté. L'application couplée d'un lemme de zéros de Philippon et d'un principe de tiroirs nous donne alors une inégalité concernant le degré d'un fermé obstructeur. 

\bigskip

Si cette inégalité ne nous permet pas encore de conclure, elle amorce la phase de descente qui constitue la cinquième partie de ce travail. Cette technique est classique dans les travaux diophantiens visant à minorer les hauteurs, mais le cas des variétés abéliennes est plus difficile, en l'absence d'un relèvement du Frobenius en caractéristique nulle. On a cherché ici à ``simuler combinatoirement'' le Frobenius, ce qui explique qu'on travaille, dès la transcendance, avec des fermés de Zariski non nécessairement irréductibles. 

\section[Préliminaires]{Préliminaires sur les variétés abéliennes}
\label{2.0}

Commençons par quelques rappels concernant la hauteur et le plongement étiré sur les variétés abéliennes. On fixe pour toute la suite une variété abélienne $A$ munie de $\mathcal{L}$ un fibré ample et symétrique, et définie sur $K$ un corps de nombres. On suppose que $A$ vérifie l'hypothèse \textbf{H} et on fixe $\mathcal{A}$ un modèle entier de $A$ sur $\mathcal{O}_K$ tel que pour une densité positive d'idéaux premiers $v$ de $\mathcal{O}_K$, la fibre spéciale $\mathcal{A}_v$ soit ordinaire.

\subsection{Choix de la constante}
\label{2.4}

La minoration que l'on a en vue fait intervenir une constante ne d\'ependant que de la variété abélienne $A$. On introduira par la suite des constantes $c_1, \ldots, c_{14}$ dépendant de $A$ et du fibré $\mathcal{L}$, mais aussi d'un certain nombre de choix sur $A$ (une base de sections globales sur $\mathcal{L}$, un modèle entier $\mathcal{A}$ donné par l'hypothèse \textbf{H}). Par abus de langage, on dira que ces constantes ne dépendent que de $A$. 

\bigskip

On choisit ensuite une constante $C_0$, ne dépendant que de $A$, grande devant les constantes $c_1, \ldots, c_{14}$ dans un sens que le cours de la preuve rendra explicitable. La constante $C(A)$ du théorème \ref{5} s'exprimera alors simplement en fonction de $C_0$:
\begin{displaymath}
C(A):=C_0^{-2 \lambda(g)},
\end{displaymath}
la dimension $g$ pouvant même être remplacée par la codimension $k$ de $V$. Toutes ces constantes sont calculables mais notre méthode ne permet pas d'améliorer la dépendance en la hauteur de $A$ du théorème \ref{4}. En effet, la constante $C_0$ est déjà au mieux monomiale en la hauteur de $A$ et le choix de la fonction $\lambda$ donne un exposant particulièrement mauvais pour la hauteur de $A$.

\subsection{Quelques définitions}
\label{2.1}

Soit $n$ un entier positif, $K$ un corps de nombres et $x:=[x_0: \cdots : x_n]$ un point de $\mathbb{P}^n(K)$. On définit la hauteur projective de $x$ par la formule suivante:
\begin{displaymath}
h(x):=\frac{1}{[K:\mathbb{Q}]} \sum_{v \in M(K)} n_v \ \mathrm{log} \big( \mathrm{max} \{ |x_0|_v, \ldots , |x_n|_v \} \big) ,
\end{displaymath} 
où $M(K)$ désigne l'ensemble des places de $K$. Les $n_v$ et les valeurs absolues sont choisis de telle sorte que la formule du produit soit satisfaite:
\begin{displaymath}
\forall y \in K, \ y \neq 0: \prod_{v \in M(K)} |y|_v^{n_v}=1,
\end{displaymath}
et que, pour tout $p$ nombre premier et $v / p$: $|p|_v=p^{-1}$. 
La formule donnant la hauteur ne dépend pas du choix d'un corps de définition $K$ de $x$, ni du choix d'un représentant de $x$ (par la formule du produit). 

\bigskip

On a fixé une variété abélienne $A$ munie de $\mathcal{L}$, un fibré ample et symétrique, et définie sur $K$, un corps de nombres. Quitte à considérer $\mathcal{L}^{\otimes 3}$, on peut supposer que $\mathcal{L}$ est projectivement normal et très ample. Il définit donc un plongement de $A$ dans un espace projectif $\mathbb{P}^n(\overline{\mathbb{Q}})$ et on fixe dès à présent une base de sections $(Y_0, \ldots, Y_n)$ de $\mathcal{L}$ sur $A$. 
Via ce plongement, pour tout point $x \in A$, on peut définir sa hauteur projective qu'on note $h_{\mathcal{L}}(x)$. 
On définit aussi sa hauteur de Néron-Tate, donnée par la formule:
\begin{displaymath}
\hat{h}_{\mathcal{L}}(x):= \lim_{n \rightarrow \infty} \frac{h_{\mathcal{L}} \left([n]x\right)}{n²},
\end{displaymath}
où $[n]$ désigne l'isogénie de multiplication par $n$, pour $n \in \mathbb{N}$. Cette fonction hauteur vérifie les propriétés:
\begin{itemize}
\item $\hat{h}_{\mathcal{L}}(x)=0$ si et seulement si $x$ est un point de torsion de $A$.
\item la hauteur est quadratique: 
\begin{displaymath}
\forall x \in A, n \in \mathbb{N}: \ \ \hat{h}_{\mathcal{L}}([n]x)= n²\hat{h}_{\mathcal{L}}(x).
\end{displaymath}
\item il existe une constante $c_{1}$ ne dépendant que de $A$ telle que: 
\begin{displaymath}
\forall x \in A, \ \ |h_{\mathcal{L}}(x)-\hat{h}_{\mathcal{L}}(x)| \leq c_{1}.
\end{displaymath}
\end{itemize}

\bigskip

Via le plongement projectif de $A$ associé à $\mathcal{L}$, on peut définir le degré relativement à $\mathcal{L}$ d'une sous-variété $V$ de $A$, qu'on note $\mathrm{deg}_{\mathcal{L}}(V)$ (voir \cite{Hartshorne77}). On a aussi défini dans l'introduction un invariant plus fin que le degré, l'indice d'obstruction. On omettra par la suite de préciser le fibré car il n'y aura pas d'ambiguïté concernant l'indice d'obstruction. 

Soit $V$ une sous-variété de $A$. La hauteur projective de $V$ relativement à $\mathcal{L}$, notée $h_{\mathcal{L}}(V)$, est définie comme la hauteur de la forme de Chow qui lui est associée via le plongement projectif. On construit ensuite la hauteur normalis\'ee $\hat{h}_{\mathcal{L}} (V)$
par passage à la limite:
\begin{displaymath}
\hat{h}_{\mathcal{L}}(V)= \mathrm{lim}_{n \rightarrow + \infty}
\frac{h_{\mathcal{L}}\left([n]V\right)\mathrm{deg}_{\mathcal{L}}(V)}{n² \ \mathrm{deg}_{\mathcal{L}}([n]V)}.
\end{displaymath}
Cette limite existe et la hauteur ainsi construite g\'en\'eralise
la hauteur de N\'eron-Tate pour les points de $A$ (voir\cite{Hindry-Silverman00}, page 450 et \cite{Philippon95}). On peut aussi définir la hauteur d´une
variété par l'approche arakelovienne (cf \cite{Bost-Gillet-Soule94}, partie 3). Ces deux notions
co\"{\i}ncident, si on fait le choix de la norme $L^2$ pour les places archimédiennes dans la définition de la hauteur projective. Le théorème des minima successifs démontré par Zhang (dans \cite{Zhang-Positive95}) implique les inégalités
suivantes:
\begin{thm}
Soit $V$ une sous-variété de $A$ de dimension $d$. Alors:
\begin{displaymath}
 \frac{1}{d+ 1} \frac{\hat{h}_{\mathcal{L}}(V)}{ \mathrm{deg}_{\mathcal{L}}(V)} \leq \me (V)
  \leq \frac{\hat{h}_{\mathcal{L}}(V)}{\mathrm{deg}_{\mathcal{L}}(V)}.
\end{displaymath}
\end{thm}

\subsection{Le plongement \'etir\'e}
\label{plongement etire}

On aura besoin de rendre négligeable la constante de comparaison entre hauteur projective et hauteur
canonique sur $A$. On résout généralement ce problème en utilisant le plongement \'etir\'e. Soit $M$ un entier strictement positif. On rappelle que l'isogénie de multiplication par $M$ sur $A$ est notée $[M]$ et on d\'efinit $\phi_M$:
\begin{eqnarray*}
  A & \rightarrow & \ \ A \times A\\
  x & \rightarrow & (x, [M] x).
\end{eqnarray*}
Cette isogénie définit un plongement de $A$ dans un espace projectif $\mathbb{P}^{m}(\overline{\mathbb{Q}})$, le {\it plongement étiré}. 

\bigskip

Le plongement étiré a été utilisé pour la première fois par Laurent dans le cadre du problème de Lehmer elliptique (\cf \cite{Laurent83}). Le principe en est le suivant: les techniques diophantiennes utilisent la hauteur projective, et le minimum essentiel fait intervenir la hauteur de Néron-Tate associée au plongement. On sait que la différence entre ces deux hauteurs est bornée mais la hauteur de Néron-Tate, qu'on prendra proche du minimum essentiel au cours de la preuve, peut être très petite. Il y a donc une perte d'information sur la hauteur projective. Le plongement étiré multiplie la hauteur par un paramètre assez grand, qui rend négligeable la constante de comparaison.

\bigskip

On note $\mathcal{M}$ le fibré associé à ce nouveau plongement et on fixe une base de sections $\textbf{Z}:=(Z_0, \ldots, Z_m)$ de $\mathcal{M}$ (en degré $1$) sur $A \times A$. On a: 
\begin{displaymath}
\mathcal{M}_{|A} \simeq \mathcal{L}^{\otimes (M^2 + 1)}. 
\end{displaymath}
Cet isomorphisme suggère le lemme suivant, qui indique la variation de la hauteur et du degr\'e par
$\phi_M$:
\begin{lem}
\label{11}
Si $x$ est un point de $A$, on a:
\begin{displaymath}
\hat{h}_{\mathcal{M}}(x)= (M^2+1)\hat{h}_{\mathcal{L}}(x).
\end{displaymath}
Si $V$ est une sous-variété de $A$, on a:
\begin{displaymath}
\mathrm{deg}_{\mathcal{M}}(V)=(M^2+1)^{\mathrm{dim}(V)}\mathrm{deg}_{\mathcal{L}}(V).
\end{displaymath}
\end{lem}
\begin{dem}
La première partie du lemme suit de la définition de $\phi_M$ par quadraticité de la hauteur.
La deuxième partie est démontrée, par exemple, dans \cite{Philippon95}, proposition 7.
\end{dem}

\subsection{Une première sélection de nombres premiers}

Pour obtenir le lemme $p$-adique crucial dans ce travail, on devra faire un certain nombre d'hypothèses sur le premier $p$ avec lequel on travaille. Les deux premières ne sont pas très restrictives et concernent les translations et dérivations sur les coordonnées projectives des variétés abéliennes. 
\begin{lem}
\label{addition}
Il existe des formes bi-homogènes de degré $(2,2)$ sur les coordonnées projectives, données par une famille finie de coefficients $\left(f_{i}\right)_{i \in I}$ telle que pour une constante $c_{2}$ ne dépendant que de $A$, pour tout $p \geq c_{2}$, tout premier $v / p$ de $\mathcal{O}_K$, on a: 
\begin{displaymath}
\forall i \in I: \ \ |f_i|_v =1.
\end{displaymath}
\end{lem}
\begin{dem}
Voir \cite{Lange-Ruppert85} ou \cite{David-Philippon02}, proposition 3.7, en remarquant qu'on a choisi $\mathcal{L}$ projectivement normal. 
\end{dem}

Un second résultat concernant les dérivations nous servira pour montrer que les coefficients du développement de Taylor d'un polynôme en les coordonnées affines sont $v$-entiers, pour les idéaux premiers $v$ considérés. 
\begin{lem}
\label{derivations algebriques}
Il existe des dérivations $(\partial_1, \ldots, \partial_g)$ sur $A$, définissant une base $(e_1, \ldots, e_g)$ du tangent en $0$, telles que pour une constante $c_{3}$ ne dépendant que de $A$, pour tout premier $p \geq c_{3}$, $v / p$ idéal premier de $\mathcal{O}_K$, on a les propriétés suivantes: 
\begin{itemize}
\item le changement de variable de la base de Siegel vers $(e_1, \ldots, e_g)$ est à coefficients dans $\mathcal{O}_K$, et l'image de $(e_1, \ldots, e_g)$ par la réduction modulo $v$ est encore une base du tangent en $0$ de $\mathcal{A}_v$;
\item pour toute base de fonctions abéliennes $(y_1, \ldots, y_n)$ définies sur un ouvert affine de $A$ (par la non-annulation d'une des sections $Y_k$ fixées précédemment), on a:  
\begin{displaymath}
\forall 1 \leq i \leq g , 1 \leq j \leq n: \ \ \ \partial_i y_j = \sum_{\alpha, \beta} a_{\alpha, \beta}^{i,j} y_{\alpha}y_{\beta}.
\end{displaymath}
avec, pour tout $(\alpha, \beta, i, j)$: $a_{\alpha, \beta}^{i,j}$ est $v$-entier.
\end{itemize}
\end{lem}
\begin{dem}
Voir \cite{Galateau08}, théorème 2.6, qui résume l'argument de \cite{David91}, théorème 4.1, où la hauteur des coefficients est bornée explicitement en fonction de $A$ dans le cas du plongement théta. Précisons que les coefficients de la matrice de l'application induite sur le tangent en $0$ par ce choix de dérivations sont dans $K$ et ont une hauteur bornée uniquement en fonction de $A$ (dans le cas du plongement théta, voir le lemme 1.4.14 de \cite{David89}). Quitte à multiplier les nouveaux paramètres par un entier assez grand (borné uniquement en fonction de $A$ et du choix de coordonnées affines), on obtient une base à coefficients dans $\mathcal{O}_K$. Le déterminant de la matrice de passage de la nouvelle base par rapport à la base de Siegel est encore de hauteur bornée par une constante ne dépendant que de $A$. 
\end{dem}
Ces deux lemmes étant posés, on peut définir l'ensemble de premiers avec lesquels on travaillera désormais:
\begin{defn}
Soit $\mathcal{P}_{\mathcal{A}}$ l'ensemble des premiers $p \geq \mathrm{max} \{ c_{2}, c_{3} \}$ tels qu'il existe un idéal premier $ v / p $ de $\mathcal{O}_K$ pour lequel $\mathcal{A}_v$ a (bonne) réduction ordinaire. Cet ensemble est de densité positive par choix du modèle.
\end{defn}

\section{Points de torsion et estimations $p$-adiques}

Dans cette partie, on donne des estimations $p$-adiques concernant certains points de torsion sur les variétés abéliennes. On commence par obtenir des propriétés concernant les paramètres en $0$ d'un point de torsion bien choisi, puis on transfère ces propriétés à des coordonnées affines (associées à un ouvert affine contenant le point de torsion) en utilisant les propriétés des dérivations algébriques. 

\subsection{Dérivations sur le plongement étiré}

On fixe ici une base de dérivations sur l'image de la variété abélienne par plongement étiré. On renvoie par exemple à \cite{David-Hindry00}, 4 pour plus de détails. 

\bigskip

Un bon choix de cette base nous permettra d'abord de transporter des estimations $p$-adiques sur un système de paramètres vers des estimations portant sur les coordonnées affines. Dans la partie suivante, on aura besoin de ce contexte différentiel pour comparer la fonction de Hilbert d'un fermé de Zariski $X$ avec multiplicité et la fonction de Hilbert d'un translaté de $X$ par un sous-groupe de torsion.  

\bigskip

On commence par choisir la base de paramètres $(t_1, \ldots, t_g)$ de $A$ en l'origine, duale de la base du tangent algébrique donnée par le lemme \ref{derivations algebriques}. A cette base de paramètres correspond une base du tangent de $A$ en $0$ qu'on note encore, par abus de langage, $(\partial_1, \ldots, \partial_g)$. Si $x$ est un point de $A$, on obtient une base de param\`etres locaux en $x$: 
\begin{displaymath}
(t_1 \circ \tau_{-x}, \ldots, t_g \circ \tau_{-x}), 
\end{displaymath}
où $\tau_{-x}$ désigne la translation par $-x$ sur $A$. Cette base de paramètres correspond à une base du tangent de $A$ en $x$, qu'on note $(\partial_{1, x}, \ldots, \partial_{g,x})$.

\bigskip

L'image de cette base par le plongement étiré fournit une base de dérivations en tout point de $\phi_M(A)$, qu'on précise ici. On note encore, par abus de langage, $(\partial_1, \ldots, \partial_g)$ l'image de la base de dérivations en l'origine de $A$ sur le premier facteur de $A \times A$, et $(\partial_{g+1}, \ldots, \partial_{2g})$ son image sur le second facteur. On a donc une base $(\partial_1, \ldots, \partial_{2g})$ de dérivations de $A \times A$ en l'origine et on en déduit par translations une base du tangent de $A \times A$ en tout point. L'application induite par $[M]$ sur le tangent de $A$ étant la multiplication par $M$, on obtient enfin une base de dérivations $(\delta_{1,x}, \ldots, \delta_{g, x})$ en tout $(x, Mx) \in \phi_M(A)$ par la formule:
\begin{displaymath}
\forall \ 1 \leq i \leq g: \ \ \delta_{i, x} = \partial_{i, (x, [M]x)} + M \partial_{g+i, (x, [M]x)}.
\end{displaymath} 
On note alors, pour $\textbf{k} = (k_1, \ldots, k_g) \in \mathbb{N}^g$:
\begin{displaymath}
|\textbf{k}|= \sum_{i=1}^g k_i \mathrm{ \ et \ } \textbf{k!}= k_1! \cdots k_g!,
\end{displaymath}
et les ``dérivées divisées'':
\begin{eqnarray*}
\mathrm{sur \ } A & : & \partial^{\textbf{k}}_x= \frac{1}{\textbf{k!}}\partial_{1, x}^{k_1} \circ \cdots \circ \partial_{g, x}^{k_g}, \\
\mathrm{sur \ } \phi_M(A) & : & \delta^{\textbf{k}}_{x} = \frac{1}{\textbf{k!}}\delta_{1, x}^{k_1} \circ \cdots \circ \delta_{g, x}^{k_g}.
\end{eqnarray*}

Si $T \geq 1$, et $\mathfrak{B}$ est un idéal homogène définissant un fermé $X$ de $A$ dans le plongement étiré, on d\'efinit $\mathfrak{B}^{(T)}$, l'id\'eal engendr\'e par les polyn\^omes $F \in \overline{\mathbb{Q}} [\textbf{Z}]$ nuls sur $X$ à un ordre au moins $T$, dans le sens suivant:
\begin{displaymath}
\forall x \in X, \ \forall \textbf{k} \in \mathbb{N}^g \mathrm{\ tq \
} |\textbf{k}| \leq T - 1: \ \ \delta^{\textbf{k}}_{x} F = 0 .
\end{displaymath}

\subsection{Coordonnées affines des points de torsion sur une variété abélienne}

Soit $p \in \mathcal{P}_{\mathcal{A}}$ et $v / p$ un idéal premier de $\mathcal{O}_K$ tel que $\mathcal{A}_v$ soit ordinaire. On a le lemme suivant: 
\begin{lem}
\label{lemme p adique}
Si $x$ est un point de $p$-torsion se réduisant sur $0$ modulo $w$, pour un idéal premier $w / v$ dans un corps de définition de $x$, on a: 
\begin{displaymath}
\forall \ 1 \leq i \leq g: \ \ \ |t_i(x)|_{w} \leq p^{-1/p}. 
\end{displaymath}
\end{lem} 
\begin{dem}
Il s'agit du Corollaire 2.8 de \cite{Galateau08}, la seule différence étant la convention faite sur les normes $w$-adiques.
\end{dem}

\bigskip

On peut maintenant démontrer la proposition qui nous sera utile dans la partie suivante. Il s'agit de passer d'une propriété sur les paramètres à une propriété concernant les sections globales de $\mathcal{M}$ sur $A$. Soit $T$ un entier, $\alpha \in A$ et $P$ un point de torsion se réduisant sur $0$ modulo $w$, où $w / v$ est une place d'un corps de définition de $P$ et de $\alpha$. On choisit un ouvert affine défini par la non-annulation d'une coordonnée projective $Y_i$ en $P$ et associé à des coordonnées affines: 
\begin{displaymath}
(y_1, \cdots , y_n).
\end{displaymath}
Soit aussi $F$ une section de $\mathcal{M}$ en degré quelconque sur $A$, à coefficients $v$-entiers sur la base de sections sur $A$ en tout degré formée par restriction à partir de $(Z_0, \ldots, Z_m)$, par projective normalité. On suppose que $F$ est nul en $\alpha$ à un ordre au moins $T$: 
\begin{displaymath}
\forall \ \textbf{k} \in \mathbb{N}^g \mathrm{ \ t.q. \ } |\textbf{k}| \leq T-1: \ \ \ \delta_{\alpha}^{\textbf{k}}F=0.
\end{displaymath}
On a alors la proposition suivante, où $\tau_{\alpha}$ désigne la translation par $\alpha$:
\begin{prop}
\label{proposition p adique}
Avec le choix de coordonnées affines induit sur le plongement étiré par $(y_1, \ldots, y_n)$, on a:
\begin{displaymath}
|F \circ \tau_{\alpha} (P)|_w \leq p^{-T/p}.
\end{displaymath}
\end{prop}
\begin{dem}
On commence par remarquer que: 
\begin{displaymath}
x \longrightarrow F (\alpha + x) 
\end{displaymath}
est une fonction rationnelle sur $A$ bien définie en $0$. En notant $\mathfrak{A}$ l'anneau des fonctions sur $A$, $\mathfrak{m}$ l'idéal maximal associé à l'origine et $\mathfrak{A}_{\mathfrak{m}}$ le localisé de $\mathfrak{A}$ en $0$, on a l'injection: 
\begin{displaymath}
\mathfrak{A}_{\mathfrak{m}} \hookrightarrow \widehat{\mathfrak{A}}_{\mathfrak{m}} \simeq K[[t_1, \ldots, t_g]],
\end{displaymath}
qui associe à une fonction régulière sur un voisinage de $0$ son développement de Taylor en les paramètres $(t_1, \ldots, t_g)$. 

\bigskip

Soit maintenant $\textbf{k} \in \mathbb{N}^g$. Commençons par démontrer que $\delta^{\textbf{k}}_{\alpha}F$, le terme d'ordre $\textbf{k}$ dans le développement en série de Taylor, est $w$-entier. Par construction, on a la formule:
\begin{displaymath}
\delta^{\textbf{k}}_{\alpha} F  = \delta^{\textbf{k}} \big(F\circ \tau_{\alpha}\big).
\end{displaymath}
Les coefficients des formes donnant l'addition sur les coordonnées projectives sont de norme $w$-adique égale à $1$ (par choix des premiers). On choisit les coordonnées projectives:
\begin{displaymath}
[1: y_1(P): \cdots: y_n(P)] 
\end{displaymath}
pour $P$ (quitte à permuter) et on peut supposer que les coordonnées projectives de $\alpha$ sont $w$-entières (quitte à multiplier par un entier assez grand). En itérant l'addition, on voit que la multiplication par n'importe quel entier est donnée sur les coordonnées projectives par des polynômes à coefficients $w$-entiers. Il suit que sur l'ouvert considéré, $F \circ \tau_{\alpha}$ est un polynôme à coefficients $w$-entiers en les coordonnées affines: 
\begin{displaymath}
\big(y_1, \ldots, y_n \big). 
\end{displaymath}
Enfin, comme les dérivations $\delta_i$ sur le plongement étiré sont des combinaisons linéaires à coefficients entiers des $\partial_i$ sur le plongement standard, il suffit de montrer que pour tout monôme unitaire $\textbf{y}^{\textbf{l}}$ ($\textbf{l} \in \mathbb{N}^n$ quelconque) en les coordonnées affines $(y_1, \ldots, y_n)$, les dérivées divisées: 
\begin{displaymath}
\partial^{\textbf{k}} \left(\textbf{y}^{\textbf{l}}\right), \ \textbf{k} \in \mathbb{N}^g,
\end{displaymath}
sont $w$-entières. Ceci découle du choix des dérivations dans le lemme \ref{derivations algebriques}. Les coefficients apparaissant quand on dérive les $y_i$ sont $w$-entiers, et comme la base de dérivations est algébrique, les factorielles dans les dérivées divisées disparaissent pour les polynômes.

\bigskip

On a donc démontré que les coefficients $\delta^{\textbf{k}}_{\alpha} F $ du développement de Taylor de $F$ en $\alpha$ sont $w$-entiers. Cette série de Taylor converge (pour la norme $w$-adique) sur la ``boule unité'': 
\begin{displaymath}
\Big\{x, \ \mathrm{sup}_{1 \leq i \leq n} \{ \ |t_i(x)|_w \} < 1 \Big\}. 
\end{displaymath}
Comme $P$ est un point de $p$-torsion se réduisant sur $0$ modulo $w$, la série de Taylor considérée converge donc en $P$ (et en particulier par le lemme \ref{lemme p adique}). On a donc, dans une complétion $w$-adique:
\begin{displaymath}
F(\alpha + P)= \sum_{\textbf{k}} \delta^{\textbf{k}}_{\alpha} \left( F \right) \cdot \textbf{t}(P)^{\textbf{k}}. 
\end{displaymath}
Puisque $F$ est nulle à l'ordre $T$ en $\alpha$, la proposition suit du lemme \ref{lemme p adique}, de la $w$-intégrité des coefficients de cette série formelle, et de l'inégalité ultramétrique. 
\end{dem}

\section{Section nulle sur une réunion de translatés de la variété}
\label{3.0}

Suivant la stratégie de \cite{Amoroso07}, on utilise la propriété $p$-adique démontrée dans le paragraphe précédent de façon ``équivariante'', dans le but de comparer deux fonctions de Hilbert, celle d'un fermé $X$ avec multiplicité $T$, et celle du translaté de $X$ par un sous-groupe de torsion sans multiplicité. Cette méthode donne directement l'existence d'une section nulle sur des translatés de $V$, alors que la précédente (voir \cite{Amoroso-David03}) commence par trouver une section nulle sur $V$ avec forte multiplicité avant d'extrapoler aux translatés. Cette approche plus ``intrinsèque'' simplifie la preuve de la façon suivante (ces trois points étant liés):  
\begin{itemize}
\item on peut se contenter de considérer les points de $p$-torsion pour un seul premier $p$; 
\item on n'est plus obligé d'itérer ``l'extrapolation'' en translatant plusieurs fois par des groupes de torsion;
\item on utilise un lemme de zéros plus simple.
\end{itemize}
 
\subsection{Le choix des paramètres}
\label{parametres}

Soit $X$ un fermé de Zariski de $A$ d'indice d'obstruction $\omega(X)$. On introduit un certain nombre de paramètres dépendant de $\omega(X)$ et de $A$. On rappelle qu'on a introduit en \ref{2.4} une constante $C_0$ ne dépendant que de $A$, grande devant toutes les constantes (ne dépendant que de $A$) du problème. Soit $\Delta$ un paramètre vérifiant l'inégalité suivante: 
\begin{displaymath}
\Delta \geq C_0 \mathrm{log} (3\omega(X)). 
\end{displaymath}
Ce paramètre permettra souvent de simplifier les calculs en éliminant les constantes dans les termes ayant une dépendance logarithmique en $\omega(X)$.

\bigskip

On choisit d'abord $p \in \mathcal{P}_{\mathcal{A}}$ ayant une dépendance d'ordre logarithmique en $\omega(X)$ dans le sens suivant: 
\begin{displaymath}
\mathrm{log}(p) \leq c_4 \mathrm{log} (\Delta),
\end{displaymath}
où on a posé: $c_{4}:= \left( 5g(g+1)\right)^{g+1}$, cette constante ne dépendant que de $A$. Soit $v / p$ un idéal de $\mathcal{O}_K$ de bonne réduction ordinaire (un tel idéal existe, par définition de $\mathcal{P}_{\mathcal{A}}$). 

\bigskip

On suppose dans toute cette partie que le minimum essentiel de $X$ est majoré de la façon suivante:
\begin{eqnarray}
\label{majoration minimum essentiel}
\me(X) <  \frac{\mathrm{log}(p) }{ \Delta p \omega(X)}.
\end{eqnarray}

\bigskip

On poursuit avec l'entier $T$ définissant la multiplicité. L'argument $p$-adique qu'on va utiliser impose que $T$ soit plus grand que $ \frac{\Delta p}{\mathrm{log}(p)}$, et l'ajustement des paramètres (pour comparer les fonctions de Hilbert, puis pour le lemme de zéros) force que $T$ soit en fait de l'ordre de ce quotient. On prend donc: 
\begin{displaymath}
T:=  \frac{ \Delta p}{\mathrm{log}(p)}.
\end{displaymath}

La comparaison des ``fonctions de Hilbert'' fait intervenir le paramètre $L$ et impose que celui-ci soit assez petit (de l'ordre d'une constante). Le choisir assez grand permet de prendre $M^2$ plus petit (puisque c'est le produit $LM^2$ qui intervient typiquement dans la preuve), donc d'améliorer la minoration du minimum essentiel. Ceci nous incite à prendre: 
\begin{displaymath}
L:= [C_0^{1/2}],
\end{displaymath}
Remarquons que le contexte du plongement étiré ``décale'' les paramètres: le choix de l'entier $L$ du cas torique correspond {\it grosso modo} à $M^2$ dans le cas abélien.

Cet entier $M$ doit être assez grand pour que la comparaison des ``fonctions de Hilbert'' considérées donne l'existence d'une section non-nulle (majoration de l'``exposant de Dirichlet'' selon la terminologie diophantienne, voir par exemple \cite{David-Hindry00}, 5.3). On doit ensuite le prendre assez petit pour obtenir une bonne majoration du degré d'un ``fermé obstructeur'' après utilisation du lemme de zéros: 
\begin{displaymath}
M:=\left[ \left(\frac{\Delta  p \omega(X)}{ \mathrm{log}(p)}\right)^{1/2}\right].
\end{displaymath}

L'inégalité (\ref{majoration minimum essentiel}) permet alors que la hauteur des petits points de $X$ dans le plongement étiré soit bornée, ce qui suffit puisqu'on travaille avec la hauteur projective, dont la différence à la hauteur de Néron-Tate est bornée par une constante. \\
\\
\textbf{Remarque} On pourrait aussi choisir une stratégie légèrement différente, et construire une section de $\mathcal{M}$ sur $A \times A$ non-nulle sur $X + \mathrm{Ker}[p]^*$ (ce qui est fait par exemple dans \cite{David-Hindry00}). Avec un autre choix de paramètres, en particulier $L$ très proche de $M^2$, on montrerait qu'une telle section est non-nulle sur l'image de $A$ par le plongement étiré, et que la dimension de l'espace des sections sur $A \times A$ est correctement minorée. Cette méthode donne ici des résultats comparables. 

\subsection{Comparaison de ``fonctions de Hilbert''}

On rappelle qu'on a fixé une base $\textbf{Z} := (Z_0, \ldots, Z_m)$ de sections de $\mathcal{M}$ sur $A \times A$ (voir \ref{plongement etire}). L'image par le plongement étiré du fermé $X$ est définie par un idéal homogène $\mathfrak{B} \subset \overline{\mathbb{Q}} [\textbf{Z}]$. On pose: 
\begin{displaymath}
\mathrm{Ker}[p]^*: = \pi_v^{-1}(0),
\end{displaymath}
où $\pi_v$ est la réduction:
\begin{displaymath}
\pi_v: \mathcal{A} \longrightarrow \mathcal{A}_v.
\end{displaymath}
On note aussi pour simplifier:  
\begin{displaymath}
X':= X+ \mathrm{Ker}[p]^*,
\end{displaymath}
et l'image de ce fermé par le plongement étiré est définie par un idéal homogène $\mathfrak{B}'$. 

\bigskip

Les morphismes de restriction à des fermés $X$ inclus dans $A$ n'étant pas nécessairement surjectifs, on est amené à poser la définition suivante:
\begin{defn}
Si $\mathfrak{I}$ est un idéal définissant un sous-schéma de $\phi_M(A)$ (non nécessairement réduit ou irréductible), on note $H(\mathfrak{I}, L)$ la dimension de l'image du morphisme de restriction:
\begin{eqnarray*}
r_{\mathfrak{I}}: H^0(A, \mathcal{M}_{|A}^{\otimes L}) & \longrightarrow & \left[\overline{\mathbb{Q}}[\textbf{Z}]/ \mathfrak{I}\right]_L \\
 f & \longrightarrow & [f ]_{\mathfrak{I}}.
\end{eqnarray*}
\end{defn}

\begin{prop}
\label{fonctions hilbert}
On a l'inégalité suivante: 
\begin{displaymath}
H(\mathfrak{B}', L) \leq 2 H(\mathfrak{B}^{(T)}, L).
\end{displaymath} 
\end{prop}
\begin{dem}
Soit $\theta$ un réel positif tel qu'on ait: 
\begin{displaymath}
\me(X) < \theta < \frac{1}{M^2};
\end{displaymath}
ceci étant possible par (\ref{majoration minimum essentiel}) et choix de $M$.
On commence par considérer l'ensemble: 
\begin{displaymath}
X(\theta) = \{ x \in X , \hat{h}_{\mathcal{L}}(x) \leq \theta \}. 
\end{displaymath}
Cet ensemble est Zariski-dense dans $X$ et l'ensemble $X'(\theta)=X(\theta)+ \mathrm{Ker}[p]^*$ est Zariski-dense dans $X'$. Pour simplifier les notations, on pose: 
\begin{eqnarray*}
h:= H(\mathfrak{B}', L) & \mathrm{et:} & h':= H(\mathfrak{B}', L)- H(\mathfrak{B}^{(T)}, L).
\end{eqnarray*}
On peut de plus supposer que $h' \geq 1$, la proposition étant déjà démontrée sinon. Pour distinguer un point $x \in A$ de son image par le plongement étiré, on notera aussi:
\begin{displaymath}
x_M:= \phi_M(x).
\end{displaymath}

\bigskip

Posons: 
\begin{displaymath}
S:= \mathrm{Ker}\left(r_{\mathfrak{B}'} \right),
\end{displaymath}
et filtrons cet espace vectoriel. Pour tout $d \in \mathbb{N}^*$, on définit: 
\begin{displaymath}
 S_d:= \Big\{ F \in S, \forall x \in X'(\theta), [K(x): K] \leq d \ : \ F(x_M)=0 \Big\}.
\end{displaymath}
Les $S_d$ forment une suite croissante de sous-espaces vectoriels de $S$. L'espace vectoriel $S$ étant de dimension finie, cette suite est stationnaire à partir d'un certain rang. On a alors, pour un entier $d_0$ assez grand: 
\begin{displaymath}
S_{d_0}= \Big\{ F \in S, \forall x \in X'(\theta) \ : \ F(x_M)=0 \Big\}.
\end{displaymath}
Comme $X'(\theta)$ est Zariski-dense dans $X'$, on a:
\begin{displaymath}
S=S_{d_0}. 
\end{displaymath}
Les points de $X'(\theta)$ étant définis sur une extension de degré au plus $d_0$ sont en nombre fini par le ``théorème de Northcott''. On note $(x_i)_{1 \leq i \leq i_0}$ ces points. On indexe aussi $(Q_j)_{1 \leq j \leq j_0}$, une base de sections de $\mathcal{M}^{\otimes L}$ sur $A$ (le plongement étant projectivement normal, il s'agit de monômes en les $Z_j$ de degré $L$). De plus, pour tout $1 \leq i \leq i_0$, on choisit $\alpha_i \in X(\theta)$ et $P_i \in \mathrm{Ker}[p]^*$ tels que: $x_i=\alpha_i + P_i$. 

\bigskip

L'orthogonal $S^{\bot}$ de $S$ dans $H^0\left(A, \mathcal{M}_{|A}^{\otimes L}\right)$, qui est de dimension $h$, est engendré par les $\big(Q_j(x_{i,M}) \big)_{i,j}$. Par définition, on peut trouver $h$ points et $h$ multi-indices tels que (quitte à renuméroter):
\begin{displaymath}
\Delta=  \mathrm{det} \big(Q_j(x_{i,M}) \big)_{1 \leq i, j \leq h}  \neq 0.
\end{displaymath}
Remarquons que ce déterminant est défini dans un espace projectif par plongements de Segré et Véronese, et que son annulation a donc un sens. On peut de plus choisir un ouvert affine pour chaque $x_i= \alpha_i + P_i$ et on choisit l'ouvert donné par la proposition \ref{proposition p adique} pour le point $P_i$ dans les calculs qui suivent. La non annulation du déterminant implique en particulier que les $(Q_j)_{1 \leq j \leq h}$ sont deux-à-deux distincts et engendrent un espace vectoriel de dimension $h$. On considère la projection: 
\begin{displaymath}
\pi:\mathrm{Vect} \big(Q_1, \dots, Q_h \big) \rightarrow \left[ \overline{\mathbb{Q}} [\textbf{Z}] / \mathfrak{B}^{(T)} \right]_{L},
\end{displaymath} 
qui n'est autre que la restriction de $r_{\mathfrak{B}^{(T)}}$ à l'espace de sections considéré. Son noyau est de dimension au moins $h'$. Ceci signifie qu'il existe $h'$ sections en degré $M$ linéairement indépendantes:
\begin{displaymath}
G_k= \sum_{j=1}^h g_{k, j} \ Q_j \ , \ k=1, \ldots, h',
\end{displaymath}
nulles sur $X$ à un ordre au moins $T$. On prend désormais $K'$ une extension finie de $K$ sur laquelle $X$ et les coordonnées des points de Ker$[p]^*$ sont définis, et $w / v$ une place ultramétrique de $K'$. Quitte à faire des opérations élémentaires sur la matrice $(g_{k, j})_{k,j}$, on peut alors supposer que:
\begin{displaymath}
G_k= \sum_{j=1}^{h-k+1} g_{k, j} \ Q_j,
\end{displaymath} 
et que, pour tout $1 \leq k \leq h'$:
\begin{eqnarray*}
|g_{k, j}|_w & \leq & 1 \ \mathrm{si} \ j < h-k+1, \\
|g_{k, j}|_w & = & 1 \ \mathrm{si} \ j=h-k+1.
\end{eqnarray*}
En utilisant ces propriétés, on remplace les $h'$ dernières colonnes de la matrice $\big(Q_j(x_{i,M})\big)_{1 \leq i, j \leq h}$ par les colonnes:
\begin{displaymath}
C_k:= \ ^t\Big( G_k(x_{1, M}), \ldots, G_k(x_{h, M}) \Big), \ \ 1 \leq k \leq h',
\end{displaymath}
en substituant d'abord $C_1$ à la dernière colonne, puis $C_2$ à l'avant-dernière, et ainsi de suite. Si on note $\Delta'$ le nouveau déterminant, on a alors:
\begin{displaymath}
|\Delta'|_w=|\Delta|_w.
\end{displaymath}
Les polynômes $G_k$ sont nuls à un ordre au moins $T$ en tous les points $\alpha_{i, M}$, pour $1 \leq i \leq h$ et ils sont à coefficients $w$-entiers.
On est donc en situation d'utiliser la proposition \ref{proposition p adique} qui donne, pour $1 \leq i \leq h$ et $1 \leq k \leq h'$: 
\begin{displaymath}
\left| G_k(x_{i, M}) \right|_w \leq p^{-\frac{T}{p}}. 
\end{displaymath}

\bigskip

On peut maintenant majorer $|\Delta'|_w$ en développant successivement le déterminant selon les $h'$ dernières colonnes. Sur la $i$-ème ligne de chacune des $h-h'$ premières colonnes, on majore la norme de chaque coefficient par: 
\begin{displaymath}
\mathrm{max} \{ 1, |z_1(x_{i, M})|_w, \ldots, |z_m(x_{i, M})|_w \}^L; 
\end{displaymath}
où les $(z_i)_{1 \leq i \leq m}$ sont un système de coordonnées affines sur un ouvert contenant $x_{i,M}$. Le système de coordonnées dépend du point, mais on omettra de préciser cette dépendance puisque la majoration fera vite apparaître la hauteur, qui ne dépend pas du choix des coordonnées affines.
On trouve finalement, par la formule donnant le déterminant en fonction de ses coefficients (sur les $h-h'$ premières colonnes) et l'inégalité ultramétrique:
\begin{displaymath}
|\Delta|_w = |\Delta'|_w \leq p^{-h'\frac{T}{p}} \prod_{i=1}^h  \mathrm{max} \left\{ 1, |z_1(x_{i, M})|_w, \ldots, |z_m(x_{i, M})|_w \right\}^{L}.
\end{displaymath} 
Si $w$ est une autre place de $K'$, on se contente d'estimations triviales et, dans le cas archimédien, la formule des déterminants fait cette fois intervenir le facteur $h!$. On a: 
\begin{displaymath}
|\Delta|_w \leq h!   \prod_{i=1}^h  \mathrm{max} \{1, |z_1(x_{i, M})|_w, \ldots, |z_m(x_{i, M})|_w \}^{L}.
\end{displaymath}

\bigskip

On écrit maintenant la formule du produit, dans sa version logarithmique, pour le déterminant $\Delta$ qui est défini sur le corps $K'$:
\begin{displaymath}
0 \leq -\frac{h' T }{[K: \mathbb{Q}]} \frac{\mathrm{log}(p)}{p} + h \mathrm{log} (h) + hL \max_{1 \leq i \leq h} \{ h_{\mathcal{M}}(x_i) \} . 
\end{displaymath}
On utilise ensuite la comparaison entre hauteur projective et hauteur de Néron-Tate (voir \ref{2.1}) et la variation de la hauteur par plongement étiré (voir le lemme  \ref{11}): 
\begin{displaymath}
\frac{h' T }{[K: \mathbb{Q}]} \frac{\mathrm{log}(p)}{p} \leq h \mathrm{log} (h) + 2hL\big(M^2\theta + c_1\big). 
\end{displaymath}
Puis, pour une constante $c_{5}$ ne dépendant que de $A$, comme $L=[C_0^{1/2}]$:
\begin{displaymath}
\frac{h'}{h} \leq c_{5} \frac{p}{T \mathrm{log}(p)} \mathrm{log}(C_0h).
\end{displaymath}
On majore le terme $\mathrm{log}(h)$ en utilisant l'inclusion triviale:
\begin{displaymath}
\mathrm{Im}(r_{\mathfrak{B}'}) \subset H^0\left(A, \mathcal{M}_{|A}^{\otimes L}\right), 
\end{displaymath}
ce qui donne:
\begin{displaymath}
\mathrm{log} ( h ) \leq c_{6} \mathrm{log} (LM^2),
\end{displaymath}
pour une constante $c_{6}$ ne dépendant que de $A$ (par le théorème de Hilbert-Serre).
On en déduit, par choix des paramètres (voir \ref{parametres}) et pour une constante $c_{7}$ ne dépendant que de $A$:
\begin{eqnarray*} 
\frac{h'}{h} & \leq & c_{7} \frac{p }{T \mathrm{log}(p)} \mathrm{log}(C_0M) \\
 & \leq & c_{7} \frac{\mathrm{log}(C_0M)}{\Delta}\\
 & \leq & c_{7} (2 + c_{4}) \frac{\mathrm{log}(\Delta)}{\Delta} + c_{7} \frac{\mathrm{log}(\omega(X))}{\Delta} \\
 & \leq & \frac{1}{2},
\end{eqnarray*}
la constante $C_0$ étant prise suffisamment grande.
La proposition est donc entièrement démontrée, par les définitions de $h$ et $h'$.
\end{dem}

\subsection{Utilisation du lemme de zéros}

On montre qu'avec le choix des paramètres, et la comparaison des fonctions de Hilbert faite précédemment, il existe une section nulle sur $X'$, le translaté de $X$ par un sous-groupe de torsion. Puis à l'aide d'un lemme de zéros classique de Philippon, on en déduit une inégalité impliquant le degré d'une variété obstructrice. Cette inégalité ne permet pas de conclure et on devra recourir à un argument de descente. A cette fin, on raffine le lemme de zéros avec un principe de tiroirs pour montrer que la variété obstructrice contient une bonne proportion de composantes du fermé algébrique d'origine, si celui-ci n'est plus irréductible.

\begin{cor}
\label{section nulle}
Il existe une section de $\mathcal{L}$ sur $A$ en degré $L(M^2+1)$, non nulle, qui s'annule sur $X'$. 
\end{cor}
\begin{dem}
On va en fait démontrer qu'il existe une section de $\mathcal{M}$ en degré $L$, non nulle sur $A$, qui s'annule sur $X'$. Ceci impliquera le corollaire par l'isomorphisme: 
\begin{displaymath}
\mathcal{M}_{|A} \simeq \mathcal{L}^{\otimes(M^2+1)}.
\end{displaymath}

On commence par majorer le membre de droite dans l'inégalité de la proposition \ref{fonctions hilbert}. 
On choisit d'abord une hypersurface $Z$ contenant $X$, réalisant $\omega(X)$ et associée à un idéal $\mathfrak{C}$. Suivant un calcul de géométrie différentielle classique (l'``astuce de Philippon-Waldschmidt'', voir par exemple \cite{Philippon-Waldschmidt88}, Lemme 6.7, \cite{David-Hindry00}, Lemme 5.1 ou \cite{Galateau08}, Proposition 5.1), on obtient: 
\begin{displaymath}
H(\mathfrak{B}^{(T)}, L) \leq H(\mathfrak{C}^{(T)}, L) \leq  T H(\mathfrak{C}, L).
\end{displaymath}
Par le théorème de Chardin et le lemme \ref{11}, on a, pour une constante $c_{8}$ ne dépendant que de $A$:
\begin{displaymath}
H(\mathfrak{C}, L) \leq \mathrm{dim}_{\overline{\mathbb{Q}}} \left[ \overline{\mathbb{Q}}[ \textbf{Z}] / \mathfrak{C} \right]_L  \leq c_{8}  (LM^2)^{g-1} \omega(X). 
\end{displaymath}
Par ailleurs, la dimension de l'espace des sections est minorée de la façon suivante (en utilisant le théorème de Hilbert-Serre), pour $c_9$ ne dépendant que de $A$:
\begin{displaymath}
\mathrm{dim}_{\overline{\mathbb{Q}}} \left( H^0(A, \mathcal{M}_{|A}^{\otimes L}) \right) \geq c_{9}  \left(LM^{2}\right)^g.
\end{displaymath}
Si toute section de $\mathcal{M}$ sur $A$ en degré $L$ nulle sur $X'$ est identiquement nulle , on a l'inégalité:
\begin{displaymath}
 c_{9} \left(LM^{2}\right)^g \leq H(\mathfrak{B}', L) \leq 2 H(\mathfrak{B}^{(T)}, L) \leq 2 c_{8}  T(LM^2)^{g-1} \omega(X). 
\end{displaymath}
Et pour une constante $c_{10}$ ne dépendant que de $A$:
\begin{displaymath}
LM^2 \leq c_{10} T \omega(X).
\end{displaymath}
Le choix des paramètres donne: 
\begin{displaymath}
\frac{1}{4}C_0^{1/2} \frac{\Delta p \omega(X)}{\mathrm{log}(p)} \leq c_{10} \frac{\Delta p \omega(X)}{\mathrm{log}(p)}, 
\end{displaymath}
ce qui est absurde. Le corollaire est donc entièrement démontré.
\end{dem}

\bigskip

Les conditions sont réunies pour écrire un lemme de zéros, et obtenir le résultat suivant concernant un ``fermé obstructeur'': 
\begin{prop}
\label{ferme obstructeur}
Soit $n_X$ le nombre de composantes irréductibles de $X$. Il existe un fermé $Z$ contenant au moins $m_X$ composantes irréductibles de $X$, avec $m_X:= [n_X/g]$, tel que:
\begin{displaymath}
\de \Big(Z+ \mathrm{Ker}[p]^* \Big)^{1/k'} \leq  \frac{\Delta^2 p \omega(X)}{\mathrm{log}(p)}.
\end{displaymath}
\end{prop}
\begin{dem}
Le corollaire \ref{section nulle} donne l'existence d'une section non nulle $F$ de $\mathcal{L}$ en degré $L(M^2+1)$ nulle sur $X+ \mathrm{Ker}[p]^*$. Considérons l'ensemble algébrique $W$ défini par les équations: 
\begin{eqnarray*}
F(x + P)=0 & , & P \in \mathrm{Ker}[p]^*.
\end{eqnarray*}
On a une décomposition de $W$ en fermés équidimensionnels: 
\begin{displaymath}
W= \bigcup_{1 \leq i \leq g} W_i,
\end{displaymath}
où $i$ désigne la codimension de $W_i$. De plus, il suit de la définition de $W$ que: 
\begin{displaymath}
\forall 1 \leq i \leq g: W_i= W_i + \mathrm{Ker}[p]^*.
\end{displaymath}
Le fermé $W$ contient $X$, donc chaque composante irréductible de $X$ est incluse dans une composante irréductible de $W$ et le principe des tiroirs de Dirichlet nous garantit qu'il existe $i_0$ tel que $W_{i_0}$ contient $m \geq \frac{n_X}{g}$ composantes irréductibles de $X$. 

\bigskip

Notons $Z$ le fermé $W_{i_0}$ et $k'$ sa codimension. Par le lemme \ref{addition}, le fermé $Z$ est incomplètement défini par des équations de degré au plus $2L(M^2+1)$ (suivant la terminologie de \cite{Philippon86}, définition 3.5, qui vaut pour des fermés éventuellement réductibles). La proposition 3.3 de \cite{Philippon86} nous donne donc, pour une constante $c_{11}$ ne dépendant que de $A$:
\begin{eqnarray*}
\de(Z)^{1/k'} & \leq & c_{11} LM^2  \\
 & \leq & c_{11} C_0^{1/2} \frac{\Delta p \omega(X)}{\mathrm{log}(p)} \\
 & \leq & \frac{\Delta^2 p \omega(X)}{\mathrm{log}(p)}.
\end{eqnarray*}
Comme $Z= Z + \mathrm{Ker}[p]^*$, la proposition est entièrement démontrée.
\end{dem}

\section{Descente finale}
\label{6.0}

Pour finir la d\'emonstration du th\'eor\`eme \ref{5}, on va d'abord préciser la proposition précédente, en choisissant un bon nombre premier $p$. Le résultat obtenu avec $X=V$ ne nous permet pas encore de conclure, car il manque une hypothèse de coprimalité entre la partie discrète du stabilisateur du fermé obstructeur (qui est une variété dans ce cas) et le premier $p$. On fait donc une descente dont le principe est le suivant: on itère la proposition \ref{ferme obstructeur} avec une suite de fermés construits à partir de $V$, et l'hypothèse de coprimalité en décalé; on obtient l'égalité des dimensions entre deux crans successifs, ce qui permet alors de ramener l'hypothèse de coprimalité au même niveau et d'obtenir enfin une contradiction.

\subsection{Choix du premier et quasi-contradiction}
\label{6.1} 
 
On va d'abord fixer le premier $p$. Celui-ci doit avoir une dépendance logarithmique en $\omega(X)$, pour que la minoration du minimum essentiel ait la précision espérée; on a d'ailleurs eu besoin de faire cette hypothèse ``raisonnable'' dans le paragraphe \ref{parametres}. Dans la descente, on construira des suites de premiers de plus en plus petits et suffisamment espacés. On commence par définir deux nouveaux paramètres utiles dans la mise en place de la descente, un réel positif $\rho$ qui permet de quantifier le choix de $p$, et un entier strictement positif $R$ qui sert à fixer l'hypothèse de coprimalité dont nous aurons besoin. On fait donc les hypothèses suivantes concernant ces deux paramètres: 
\begin{displaymath}
\Delta \geq \mathrm{max} \left\{C_0 \mathrm{log}\big((3 \de(X) \big), \mathrm{log} (R) \right\} \ \  \mathrm{et:} \ \ 2 \leq \rho \leq \big(5g(k+1)\big)^{k+1}.
\end{displaymath}

\bigskip

Soit $n_X$ le nombre de composantes irréductibles de $X$ et $m_X:=\left[n_x/g \right]$. On obtient la ``quasi contradiction'':
\begin{prop}
\label{29}
On suppose que $X$ n'est pas incluse dans le translaté d'une sous-variété abélienne et que son minimum essentiel est majoré de la façon suivante:
\begin{displaymath}
\me(X) \ \omega(X) < \frac{1}{\Delta^{ \rho + 1}}.
\end{displaymath}
Alors il existe un premier $p$ strictement positif et un fermé de Zariski $Z$, de codimension $k'$, contenant au moins $m_X$ composantes de $X$, tels que:
\begin{itemize}
\item L'entier $p$ est premier avec $R$ et:
\begin{displaymath}
\frac{\Delta^{\rho}}{2}  \leq p \leq \Delta^{\rho}.
\end{displaymath}
\item On a l'inégalité: 
\begin{displaymath}
\de \Big(Z+ \mathrm{Ker}[p]^*\Big)^{1/k'} < \Delta^{2} p \omega(X). 
\end{displaymath}
\end{itemize}
\end{prop}
\begin{dem}
Il s'agit essentiellement de démontrer qu'il existe un premier $p \in \mathcal{P}_{\mathcal{A}}$ premier à $R$ et vérifiant les inégalités annoncées. Notons $n_{\Delta}$ le nombre de premiers $p$de $\mathcal{P}_{\mathcal{A}}$ vérifiant: 
\begin{displaymath}
\Delta^{\rho}/2 \leq p \leq \Delta^{\rho}, 
\end{displaymath}
et $n_R$ le nombre de premiers de $\mathbb{Z}$ divisant $R$. Les premiers de $\mathcal{P}_{\mathcal{A}}$ sont en densité $c_{12} >0$ fixée avec $A$, et comme $\Delta \geq C_0$, quitte à prendre cette constante assez grande, le théorème des nombres premiers nous donne: 
\begin{displaymath}
n_{\Delta} \geq c_{12} \frac{\Delta^{\rho}}{3\rho \mathrm{log}(\Delta)}; 
\end{displaymath}
de plus, on a: 
\begin{displaymath}
R= \prod_{i=1}^{n_R} p_i^{\nu_i} \geq \prod_{i=1}^{n_R} p_i \geq 2^{n_R}.
\end{displaymath}
On trouve donc: 
\begin{eqnarray*}
n_{\Delta} -n_R & \geq & c_{12} \frac{\Delta^{\rho}}{3\rho \mathrm{log}(\Delta)} - \frac{\mathrm{log}(R)}{\mathrm{log}(2)} \\
 & \geq & c_{12} \Delta \left( \frac{\Delta^{\rho-1}}{3\rho \mathrm{log}(\Delta)} - \frac{1}{\mathrm{log}(2)} \right) \\
 & \geq & c_{12} \Delta;
\end{eqnarray*}
car $\rho \geq 2$, et avec $C_0$ assez grande. Pour la même raison et comme $c_{12} > 0$, on a finalement: 
\begin{displaymath}
n_{\Delta} - n_R \geq 1.
\end{displaymath}
L'existence d'un premier $p \in \mathcal{P}_A$ vérifiant l'encadrement prescrit et ne divisant pas $R$ est donc acquise. 

\bigskip

Pour appliquer la proposition précédente, il reste à voir que l'hypothèse faite ici sur le minimum essentiel est plus restrictive que (\ref{majoration minimum essentiel}). La majoration de $p$ donne: 
\begin{displaymath}
\me(X) \ \omega(X) < \frac{1}{\Delta^{ \rho + 1}} \leq \frac{1}{\Delta p} \leq \frac{\mathrm{log}(p)}{\Delta p}.
\end{displaymath}
On peut donc appliquer la proposition \ref{ferme obstructeur}, qui donne l'existence du fermé $Z$ vérifiant les propriétés annoncées.
\end{dem}
\textbf{Remarque} Expliquons pourquoi cette proposition donne presque la contradiction attendue. Posons $X=V$ et supposons de plus que $V$ n'est pas inclus dans le translaté d'une sous-variété abélienne propre de $A$. Le fermé $Z$ donné par la proposition est alors irréductible. Supposons qu'on puisse imposer:
\begin{displaymath}
p \nmid [\mathrm{Stab}(Z): \mathrm{Stab}(Z)^0], 
\end{displaymath}
où $\mathrm{Stab}(Z)$ désigne le stabilisateur de $Z$, et $\mathrm{Stab}(Z)^0$ la composante connexe de $\mathrm{Stab}(Z)$ contenant l'origine. On a: 
\begin{displaymath}
\de \big(Z+ \mathrm{Ker}[p]^*\big) = \frac{p^g}{|\mathrm{Ker}[p]^* \cap \mathrm{Stab}(Z)|}\de(Z), 
\end{displaymath}
et seule la composante neutre du stabilisateur interviendrait au dénominateur. De plus, la variété $Z$, contenant $V$, n'est pas incluse dans un translaté de sous-variété abélienne, ce qui entraîne: 
\begin{displaymath}
|\mathrm{Ker}[p]^* \cap \mathrm{Stab}(Z)| \leq p^{\mathrm{dim} \left( \mathrm{Stab} (Z)^0 \right)} \leq p^{g-k'-1}
\end{displaymath}
On aurait finalement:
\begin{displaymath}
\Delta^{\rho / k'} \big( \de (Z) \big)^{1/k'} < 2 \Delta^2 \omega(V),
\end{displaymath}
Une contradiction suivrait immédiatement en prenant $\rho=3k$, et en utilisant le théorème de Chardin (voir \cite{Chardin89}, 2, Corollaire 2) pour extraire une hypersurface $\tilde{Z}$ contenant $Z$ de degré:
\begin{displaymath}
\de(\tilde{Z}) \leq c_{13} \de(Z)^{1/k'}, 
\end{displaymath}
pour $c_{13}$ ne dépendant que de $A$. On obtiendrait:
\begin{displaymath}
\me(V) \geq \frac{C(A)}{\omega(V)} \mathrm{log}\big(3 \de(V)\big)^{-4k} ,
\end{displaymath}
pour une constante $C(A)$ ne dépendant que de $A$.

\subsection{Descente en codimension générale}

La construction pr\'ec\'edente ne  permet pas de conclure, et on est amen\'e \`a it\'erer, au plus $k$ fois, la derni\`ere proposition. Cet argument de descente a déja été utilisé par Amoroso et David dans \cite{Amoroso-David99} et dans \cite{Amoroso-David03}. Il est ici compliqué par l'absence d'isogénies relevant le Frobenius en caractéristique nulle, qui rend la preuve bien plus simple dans le cas torique. Son absence pousse à le ``reconstituer combinatoirement'' en travaillant avec des réunions de tranlatés, et ceci explique qu'on ait écrit la partie diophantienne avec des fermés de Zariski au lieu de variétés. 

\bigskip
 
On va donc itérer la construction du fermé $Z$ donné par la proposition \ref{29} (et en déduire l'existence d'une bonne variété), en faisant à chaque étape une hypothèse de coprimalité entre l'entier $p$ et le stabilisateur du fermé construit au cran précédent. On pourra de plus imposer un ``emboîtement'' entre ces variétés et on démontrera qu'il existe une suite de variétés ainsi construites associée à une suite de dimensions qui n'est pas strictement croissante. La comparaison des degrés, avec l'hypothèse de coprimalité, permettra cette fois-ci de conclure. 

\bigskip

Commençons par préciser les paramètres propres à la descente. Posons:
\begin{displaymath}
\forall 1 \leq i \leq k \ , \ \rho_i:=\big(5g(k+1)\big)^{k+1-i};
\end{displaymath}
les $\rho_i$ sont décroissants avec $i$. Notons aussi:
\begin{displaymath}
P_i:= \Delta^{\rho_i},
\end{displaymath}
où: 
\begin{displaymath}
\Delta(V)= C_{0}^2 \mathrm{log}((3 \de(V)).
\end{displaymath}

\bigskip

On fixe une variété $V$ qui n'est pas incluse dans un translaté de sous-variété abélienne propre de $A$.
On va itérer la proposition \ref{29} au plus $k$ fois, en employant le paramètre $\rho_i$ à la $i$-ème étape. Le paramètre $P_i$ permettra d'encadrer le premier $p_i$ construit. On remarque que le choix des $\rho_i$ et de $\Delta$ est cohérent avec les hypothèses faites au début de \ref{6.1}. 

\bigskip

Il nous faut un nouvel indice d'obstruction pour traiter la perte de composantes dans le lemme de zéros (dont on a gardé une proportion d'au moins $\frac{1}{g}$ par un principe de tiroirs). Si $\Sigma$ est un sous-ensemble fini de $A$ et $V$ une sous-variété de $A$, on commence donc par noter:
\begin{displaymath}
\alpha_V(\Sigma):= \mathrm{nombre \ de \ composantes \ distinctes \ de \ V + \Sigma}.
\end{displaymath}

\begin{defn}
On définit l'indice d'obstruction de $V$ relativement à un sous-ensemble fini $\Sigma$ de $A$:
\begin{displaymath}
\omega(V, \Sigma):= \inf_{\begin{array}{rcl} \Lambda \subset \Sigma, \Lambda \neq \emptyset \\ V + \Lambda  \subset  Z \\ \end{array}} \left\{ \frac{\alpha_V(\Sigma)}{\alpha_V(\Lambda)} \de (Z) \right\},
\end{displaymath}
où $Z$ n'est pas nécessairement irréductible.
\end{defn}  
Cette définition permet de prendre en compte la perte de translatés de $V$ en compensant cette perte par une chute du degré de la variété obstructrice correspondante.

\bigskip

On construit des suites de vari\'et\'es ayant les propri\'et\'es suivantes:
\begin{defn}
\label{def 7}
Soit $\mathcal{W}(V)$ l'ensemble des quadruplets $(s, \textbf{p}, {\bf \Sigma }, \textbf{W})$ o\`u $s$ est un entier tel que $1 \leq s \leq k$, $\textbf{p}=(p_1, \ldots p_s)$ est un $s$-uplet de nombres premiers tels que: 
\begin{eqnarray*}
\forall \ 1 \leq i \leq s & : & \frac{P_i}{2} \leq p_i \leq P_i, 
\end{eqnarray*}
l'ensemble ${ \bf \Sigma}= (\Sigma^1, \ldots, \Sigma^s)$ est un $s$-uplet de sous-ensembles finis de $A$ et $\textbf{W}=(W_0, \ldots, W_s)$ est un $(s+1)$-uplet de sous-vari\'et\'es propres et irr\'eductibles de $A$ avec $V = W_0$, et les propriétés:
\begin{itemize}
\item pour tout $1 \leq i \leq s$, l'entier $p_i$ est premier avec $[\mathrm{Stab}(W_{i-1}): \mathrm{Stab}(W_{i-1})^0]$; $\Sigma^i$ est un ensemble non vide tel que: 
\begin{displaymath}
W_{i-1} + \Sigma^i \subset W_i;
\end{displaymath}
\item on pose $V_0=V$ et pour tout $1 \leq i \leq s$, on définit: 
\begin{displaymath}
V_i:= V + \Sigma^1 + \cdots + \Sigma^i \ \mathrm{et} \ \Lambda^i:= \mathrm{Ker}[p_i]^*.
\end{displaymath}
Alors on a:
\begin{displaymath}
\Big( \frac{\alpha_{W_{i-1}}(\Lambda^i)}{\alpha_{W_{i-1}}(\Sigma^i)} \de(W_i)\Big)^{1/\mathrm{codim}(W_i)} \leq p_i \omega( V_{i-1}) \Big( \Delta^6 (P_{i+1} \cdots P_k)^{g} \Big);
\end{displaymath}
\item on a la ``relation de récurrence'' entre indices d'obstruction, pour tout $1 \leq i \leq s$ et $\Gamma^{i} \subset \Lambda^{i}$:
\begin{displaymath}
\omega(V, \Sigma^1 + \cdots + \Sigma^{i-1} + \Gamma^{i}) \leq \Delta^3 p_{i} \omega(V_{i-1}).
\end{displaymath}
\end{itemize}
\end{defn}
\textbf{Remarque} Au sujet du dernier point, on convient que si $i=1$: 
\begin{displaymath}
\Sigma^1 + \cdots + \Sigma^{i-1} + \Gamma^i= \Gamma^i.
\end{displaymath}

\bigskip

On utilisera le lemme suivant, qui compare grossièrement deux indices d'obstruction successifs:
\begin{lem}
\label{30}
Si $(s, \textbf{p}, {\bf \Sigma }, \textbf{W}) \in \mathcal{W}(V)$, pour tout $1 \leq i \leq s$:
\begin{displaymath}
\omega(V_i) \leq p_i^g \omega(V_{i-1}).
\end{displaymath}
\end{lem}
\begin{dem}
Soit $Z$ une hypersurface, pas nécessairement irréductible, contenant $V_{i-1}$ et réalisant $\omega(V_{i-1})$. On a:
\begin{displaymath}
V_i \subset Z + \Lambda^i;
\end{displaymath}
et le lemme suit d'une majoration grossière de deg$(Z + \Lambda^i)$.
\end{dem}

\noindent
Notons:
\begin{displaymath}
\mathcal{W}_0 (V)= \{ (s, \textbf{p}, {\bf \Sigma}, \textbf{W}) \in \mathcal{W}(V), \mathrm{dim}(W_0) < \cdots < \mathrm{dim}(W_s) \}.
\end{displaymath}
Si on trouve un quadruplet $(s, \textbf{p}, {\bf \Sigma}, \textbf{W}) \in \mathcal{W} \setminus \mathcal{W}_0$, la comparaison des degr\'es aux vari\'et\'es de même dimension permet de trouver une contradiction:
\begin{prop}
\label{31}
On a l'égalité: $\mathcal{W}_0(V) = \mathcal{W}(V)$.
\end{prop}
\begin{dem}
Supposons par l'absurde qu'il existe $(s,\textbf{p}, {\bf \Sigma}, \textbf{W}) \in \mathcal{W}(V) \setminus \mathcal{W}_0(V)$. On peut donc trouver $1 \leq i \leq s$ tel que:
\begin{displaymath}
\mathrm{dim} (W_{i-1}) = \mathrm{dim} (W_i) \ \ \mathrm{et:}  \ W_{i-1} + \Sigma^i \subset W_i.
\end{displaymath} 
On observe d'abord que pour tout $0 \leq j \leq s$, $W_j$ contient un translaté de $V$. Comme $V$ n'est inclus dans aucun translat\'e d'une sous-vari\'et\'e ab\'elienne propre, il en va donc de même pour $W_j$. 

\bigskip

Par irréductibilité de $W_i$, on observe que $\alpha_{W_{i-1}}(\Sigma^i)=1$. De plus, l'entier $p_i$ est premier avec $[\mathrm{Stab}(W_{i-1}): \mathrm{Stab}(W_{i-1})^0]$, la variété $W_{i-1}$ n'est pas un translaté de sous-variété abélienne, et on a:
\begin{displaymath}
\alpha_{W_{i-1}}(\Lambda_i) \geq p_i^{1+\mathrm{codim} (W_{i})}.
\end{displaymath}
On en d\'eduit:
\begin{displaymath}
p_i \Big( p_i \de(W_{i-1})\Big)^{1/\mathrm{codim}(W_{i})} \leq p_i \omega(V_{i-1}) \Delta^6 (P_{i+1} \cdots P_k)^{g}.
\end{displaymath}
Mais par choix des paramètres: 
\begin{displaymath}
p_i^{-1/\mathrm{codim}(W_i)} \Delta^6 (P_{i+1} \ldots P_{k})^{g} \leq 2 \Delta^{u},
\end{displaymath}
avec, en remarquant que les $\rho_i$ sont une suite géométrique et $g \geq 2$:
\begin{eqnarray*}
ku & = & -\rho_i + 6k + kg(\rho_{i+1}+ \cdots + \rho_k) \\
 & \leq & -\rho_i + 5k g \rho_{i+1} \\
 & \leq & - 5(k+1)g \rho_{i+1} + 5 kg \rho_{i+1}<0.
\end{eqnarray*}
Une récurrence immédiate montre aussi que: 
\begin{displaymath}
V_{i-1} \subset W_{i-1}.
\end{displaymath} 
Puisque l'indice d'obstruction avec poids est non nul, on en déduit une contradiction.
\end{dem}

\subsection{Un dernier argument g\'eom\'etrique}

Dans le paragraphe précédent, on avait seulement besoin de savoir que $V$ n'est pas incluse dans un translaté de sous-variété abélienne. On rajoute maintenant l'hypothèse suivante, concernant son minimum essentiel:
\begin{displaymath}
\me(V) \ \omega(V) < \frac{1}{\Delta^{\left(5g(k+1)\right)^{k+1}}}.
\end{displaymath}
La preuve du théorème \ref{5} sera quasiment achev\'ee si on r\'eussit \`a d\'emontrer la négation de la proposition \ref{31}, \`a savoir:
\begin{prop}
\label{32}
On a: $\mathcal{W}_0(V) \neq \mathcal{W}(V)$.
\end{prop}
Pour y arriver, on a besoin de quelques préliminaires. A un \'el\'ement $(s,\textbf{p}, {\bf \Sigma}, \textbf{W})$ de $\mathcal{W}$, on peut associer la suite de longueur $s+1$, croissante, des dimensions des $W_i$.
On commence par construire un ordre total sur ces suites. Soient $v=(v_i)_{0 \leq i \leq s}$ et $v'=(v'_{i})_{0 \leq i \leq s'}$ deux suites croissantes d'entiers positifs.
\begin{defn}
On dit que $v \preceq v'$ si $ (v_i)_{0 \leq i \leq \mathrm{min} \{s,s'\}} < (v'_{i})_{0 \leq i \leq \mathrm{min} \{s,s'\}}$ pour l'ordre lexicographique ou, en cas d'\'egalit\'e, si $s \geq s'$.
\end{defn}
\textbf{Remarque} Si une suite est la tronquée d'une autre suite, elle est plus grande selon cet ordre (il y a un renversement des inégalités).

\bigskip

On va d\'emontrer la proposition \ref{32} par un argument de descente. On prendra un \'el\'ement $(s, \textbf{p}, {\bf \Sigma}, \textbf{W}) \in \mathcal{W}_0(V)$ associ\'e \`a une suite minimale, et on construira un nouvel \'el\'ement $(s', \textbf{p}', {\bf \Sigma}', \textbf{W}')$ dans  $\mathcal{W}(V)$, associ\'e \`a une suite strictement inf\'erieure. On sera alors assur\'e que $(s', \textbf{p}', {\bf \Sigma}', \textbf{W}')$ est dans $\mathcal{W}(V) \setminus \mathcal{W}_0(V)$.

\bigskip

Comme il n'y a qu'un nombre fini de suites de $k+1$ entiers compris entre $0$ et $k$, les suites croissantes de dimensions associ\'ees aux \'el\'ements de $\mathcal{W}_0(V)$ sont en nombre fini. De plus: 
\begin{lem}
\label{non vide}
On a:
\begin{displaymath}
\mathcal{W}_0(V) \neq \emptyset. 
\end{displaymath}
\end{lem}
\begin{dem}
Compte tenu de la proposition \ref{31}, on peut se contenter de trouver un élément dans $\mathcal{W}(V)$. On remarque d'abord (voir \cite{Amoroso-David99}, 2 qui se transpose littéralement au cas abélien): 
\begin{displaymath}
\mathrm{log}  \left( [\mathrm{Stab}(V):\mathrm{Stab}(V)^0] \right) \leq \mathrm{log} \Big( \de (\mathrm{Stab}(V))\Big) \leq g \mathrm{log} \Big(\de(V)\Big) \leq \Delta.
\end{displaymath}

\bigskip

La procédure diophantienne (proposition \ref{29}) avec $X=V$ et $R=[\mathrm{Stab}(V):\mathrm{Stab}(V)^0]$ donne l'existence d'un entier $p_1$ tel que: 
\begin{displaymath}
\frac{P_1}{2} \leq p_1 \leq P_1,
\end{displaymath}
 et d'un fermé $W_1$ contenant $W_0 + \Lambda^1$ tel que:
\begin{displaymath}
\de (W_1) ^{1/\mathrm{codim}(W_1)} \leq  p_1 \Delta^2 \omega(V).
\end{displaymath}
De plus, le premier $p_1$ est premier à $[\mathrm{Stab}(W_0):\mathrm{Stab}(W_0)^0]$ par construction. 

\bigskip

La relation de récurrence entre indices d'obstruction, pour tout $\Gamma^1 \subset \Lambda^1$, vient en extrayant une hypersurface $Z_1$ contenant $W_1$ par le théorème de Chardin (voir \cite{Chardin89}, 2, Corollaire 2): 
\begin{displaymath}
\omega(V, \Gamma^1) \leq \de(Z_1) \leq p_1 \Delta^3 \omega(V).
\end{displaymath}
La constante de comparaison entre les degrés dans le théorème de Chardin est éliminée en remplaçant $\Delta^2$ par $\Delta^3$.
\end{dem}

On peut donc choisir une suite minimale pour l'ordre $\preceq$, associ\'ee \`a un quadruplet $(s, \textbf{p}, {\bf \Sigma}, \textbf{W})$.  Le lemme technique suivant montre comment construire l'\'el\'ement $(s', \textbf{p}', {\bf \Sigma}', \textbf{W}')$ en utilisant le théorème de Bézout. Il restera ensuite à prouver que ce quadruplet est encore dans $\mathcal{W}(V)$.

\begin{lem}
\label{33}
Il existe un entier $\frac{P_{s+1}}{2} \leq p_{s+1} \leq P_{s+1}$ premier avec $[\mathrm{Stab}(W_s):\mathrm{Stab}(W_s)^0]$, un indice $1 \leq s' \leq s+1$, un sous-ensemble $\tilde{\Sigma}^{s'}$ de $\Lambda^{s'}$ et une sous-vari\'et\'e $Z_{s'}$ propre et irr\'eductible, dont le degr\'e vérifie:
\begin{displaymath}
\frac{\alpha_{W_{s'-1}}(\Lambda^{s'})}{\alpha_{W_{s'-1}}(\tilde{\Sigma}^{s'})} \de(Z_{s'}) \leq  \frac{\alpha_{W_{s'-1}}(\Lambda^{s'})}{\alpha_{W_{s'-1}}(\Sigma^{s'})} \de(W_{s'}) p_{s'} \omega( V_{s'-1}) \Big(\Delta^6 p_{s'+1}^g \cdots p_{s+1}^g \Big) ,
\end{displaymath}
et telle qu'on ait l'inclusion:
\begin{displaymath}
W_{s'-1} + \tilde{\Sigma}^{s'} \subset Z_{s'}, 
\end{displaymath} 
avec: 
\begin{displaymath}
\mathrm{codim}(Z_{s'})= \mathrm{codim}(W_{s'})+1, 
\end{displaymath}
en posant les conventions: $\mathrm{codim}(W_{s+1})=0$, $\de(W_{s+1})=1$ et $\Sigma^{s+1}= \Lambda^{s+1}$.
\end{lem}
\noindent
\textbf{Remarque}
Avec la définition de l'ordre $\prec$, on vient bien de construire une suite de dimensions ``plus petite'':
\begin{displaymath}
\Big( \mathrm{dim}(W_0), \ldots, \mathrm{dim}(W_{s'-1}), \mathrm{dim}(Z_{s'})\Big) \prec \Big( \mathrm{dim}(W_0), \ldots, \mathrm{dim}(W_s) \Big).
\end{displaymath}
\begin{dem}
On va utiliser la proposition \ref{29} avec:
\begin{displaymath}
\rho=\rho_{s+1} \ \mathrm{et} \  R=[\mathrm{Stab}(W_s):\mathrm{Stab}(W_s)^0] \prod_{i=1}^s p_i;
\end{displaymath}
on pose:
\begin{displaymath}
\Sigma_0= \Sigma^1 + \cdots + \Sigma^{s},
\end{displaymath}
puis on prend $\Sigma$ tel qu'il existe $Z$ contenant $V + \Sigma$ et vérifiant:
\begin{displaymath}
\omega(V, \Sigma_0)= \frac{\alpha_V(\Sigma_0)}{\alpha_V(\Sigma)} \de(Z).
\end{displaymath}
On choisit alors: 
\begin{displaymath}
X= V + \Sigma.
\end{displaymath}

Le fermé $X$ étant réunion d'un nombre fini de translatés de $V$ par des points de hauteur nulle, on a:
\begin{displaymath}
\me(X) = \me(V).
\end{displaymath}
De plus, le lemme \ref{30} it\'er\'e donne:
\begin{displaymath}
\omega(X) \leq \omega(V_s) \leq p_1^g \cdots p_s^g \omega(V) \leq P_1^g \cdots P_s^g \ \omega(V).
\end{displaymath}
Il en r\'esulte:
\begin{eqnarray*}
\omega(X) \me(X) & \leq & P_1^{g} \cdots P_s^{g} \omega(V) \me(V) \\
 & \leq & \Delta^{u}, \\
\end{eqnarray*}
où: 
\begin{eqnarray*}
u & = & -\left( 5g(k+1)\right)^{k+1} +g(\rho_1 + \cdots + \rho_s)\\
 & \leq & -\left( 5g(k+1)\right)^{k+1} +g(\rho_1 + \cdots + \rho_{s+1}) - \rho_{s+1}\\
 & \leq & \rho_1\left( 2g-5g(k+1)\right) - \rho_{s+1} \\
 & < & -1 - \rho_{s+1}. \\
\end{eqnarray*}
L'hypothèse de la proposition \ref{29} est donc v\'erifi\'ee avec $\rho= \rho_{s+1}$. De plus:
\begin{displaymath}
C_0 \mathrm{log} (3 \omega(X)) \leq C_0 \Big( g \mathrm{log} P_1 + \cdots + g \mathrm{log} P_s + \mathrm{log}(3 \omega(V))\Big) \leq \Delta(V).
\end{displaymath}
On majore ensuite log$(R)$, en remarquant que: 
\begin{displaymath}
\mathrm{log} (R) \leq \mathrm{log \ deg} \big( \mathrm{Stab} (W_s) \big) + \mathrm{log} (P_1) + \cdots + \mathrm{log} (P_s),
\end{displaymath}
puis de la même façon que dans la preuve du lemme \ref{non vide}:
\begin{eqnarray*}
\mathrm{log} (R) & \leq & (\mathrm{dim}(W_s)+1) \mathrm{log} \big(\de(W_s)\big) + \frac{1}{2} \Delta(V) \\
& \leq & g \Big(g \big(\mathrm{log} \Delta+ \mathrm{log} P_{s} + \cdots + \mathrm{log} P_k \big)+ g \mathrm{log} \big( \omega(V_{s-1}) \big) \Big) + \frac{1}{2} \Delta(V) \\
& \leq & \Delta(V).
\end{eqnarray*}

\bigskip

La proposition \ref{29} donne donc l'existence d'un entier: 
\begin{displaymath}
\frac{P_{s+1}}{2} \leq p_{s+1} \leq P_{s+1}, 
\end{displaymath}
strictement positif, d'un fermé de Zariski $Z$, de codimension $k'$ et contenant au moins $[\frac{\alpha_V(\Sigma)}{g}]$ composantes de $X$, tel que: 
\begin{displaymath}
\de \left(Z + \Lambda^{s+1}\right)^{1/k'} < \Delta^2 p_{s+1} \omega(X). 
\end{displaymath}
De plus, le résultat d'interpolation de Chardin (voir \cite{Chardin89}, 2, Corollaire 2) permet d'extraire une hypersurface non nécessairement irréductible $\tilde{Z}$ contenant $Z + \Lambda^{s+1}$ telle que, pour une constante $c_{14}$ ne dépendant que de $A$:
\begin{eqnarray}
\label{hypersurface extraite}
\de(\tilde{Z}) \leq c_{14} \de \left( Z + \Lambda^{s+1} \right)^{1/k'} \leq c_{14} \Delta^2 p_{s+1} \omega(X). 
\end{eqnarray}

\bigskip
 
On peut maintenant démontrer l'inégalité entre les indices d'obstruction au ``cran $s+1$''. Remarquons que si $\Gamma^{s+1}$ est un sous-ensemble quelconque de $\Lambda^{s+1}$, on a:
\begin{eqnarray*}
\omega(V, \Sigma_0 + \Gamma^{s+1}) & \leq & 2g \frac{\alpha_V(\Sigma_0+ \Gamma^{s+1})}{\alpha_V(\Sigma + \Gamma^{s+1})} \de(\tilde{Z}) \\
 & \leq & \Delta^3 p_{s+1} \frac{\alpha_V(\Sigma_0)\alpha_V(\Gamma^{s+1})}{\alpha_V(\Sigma)\alpha_V(\Gamma^{s+1})} \omega(V + \Sigma). \\
\end{eqnarray*}
On a fait apparaître les facteurs $\alpha_V(\Gamma^{s+1})$ en remarquant que les points de $\Gamma^{s+1}$ sont d'ordre premier aux ordres des points de $\Sigma^i$, pour tout $1 \leq i \leq s$. On utilise alors les propriétés fines des nouveaux indices d'obstruction, et la définition de $\Sigma$, qui garantit l'existence d'une hypersurface $Z$ contenant $V + \Sigma$ telle que: 
\begin{displaymath}
\omega(V + \Sigma_0) \geq \omega(V, \Sigma_0)=\frac{\alpha_V(\Sigma_0)}{\alpha_V(\Sigma)} \de(Z) \geq \frac{\alpha_V(\Sigma_0)}{\alpha_V(\Sigma)} \omega(V + \Sigma).
\end{displaymath}
On en déduit l'inégalité suivante:
\begin{eqnarray}
\label{recurrence obstruction}
\omega(V, \Sigma_0 + \Gamma^{s+1}) & \leq & \Delta^3 p_{s+1} \omega(V+ \Sigma_0).
\end{eqnarray}

\bigskip

{\it Premier cas.} L'hypersurface $\tilde{Z}$ contient $W_s + \Lambda^{s+1}$. On prend une composante irréductible $Z_{s+1}$ de $\tilde{Z}$ ``au-dessus de la moyenne'' dans le sens suivant. Il existe $\tilde{\Sigma}^{s+1}$ un sous-ensemble fini et non-vide de $\Lambda^{s+1}$ tel qu'on a: 
\begin{displaymath}
W_{s} + \tilde{\Sigma}^{s+1} \subset Z_{s+1},
\end{displaymath}
et que:
\begin{eqnarray}
\label{bonne composante}
\frac{\de(Z_{s+1})}{\de(\tilde{Z})} \leq \frac{\alpha_{W_s}(\tilde{\Sigma}^{s+1})}{\alpha_{W_s}(\Sigma^{s+1})}.
\end{eqnarray}

Si on note $\tilde{Z}= \bigcup_{\beta} Z_{\beta}$ la décomposition de $\tilde{Z}$ en composantes irréductibles et si $\Sigma^{(\beta)}$ est l'ensemble des $x \in \Sigma^{s+1}$ tels que $x+ W_{s} \subset Z_{\beta}$, on a en effet:
\begin{displaymath}
\frac{\sum_{\beta} \de(Z_{\beta})}{\de(\tilde{Z})} \ = \ \frac{\sum_{\beta} \alpha_{W_s} (\Sigma^{(\beta)})}{\alpha_{W_s}(\Sigma^{s+1})} \ = \ 1,
\end{displaymath}
et l'existence de $Z_{s+1}$ vérifiant (\ref{bonne composante}) suit. L'inégalité du lemme est alors démontrée avec $s'=s+1$, en remarquant que $\Lambda^{s+1}= \Sigma^{s+1}$:
\begin{displaymath}
\frac{\alpha_{W_s}(\Sigma^{s+1})}{\alpha_{W_s}(\tilde{\Sigma}^{s+1})} \de(Z_{s+1}) \leq \Delta^3 p_{s+1} \omega(X) \leq \Delta^3 p_{s+1} \omega(V_s).
\end{displaymath}

\bigskip

{\it Deuxième cas.} Le fermé $\tilde{Z}$ ne contient pas: 
\begin{displaymath}
W_s + \Lambda^{s+1}= W_s + \Sigma^{s+1}. 
\end{displaymath}
Remarquons que translater $\tilde{Z}$ ne change rien à son degré. Soit $s'$ le plus grand entier $\leq s$ tel qu'il existe un ensemble $\Gamma^{s'}$ avec (quitte à translater $\tilde{Z}$ par un point de torsion): 
\begin{displaymath}
W_{s'-1} + \Gamma^{s'} \subset \tilde{Z},
\end{displaymath}
vérifiant:
\begin{eqnarray}
\label{clause technique}
\frac{\alpha_{W_{s'-1}}(\Gamma^{s'})}{\alpha_{W_{s'-1}}(\Sigma^{s'})} \geq \frac{1}{2g} \frac{\alpha_V(\Sigma)}{\alpha_V({\Sigma_0})};
\end{eqnarray}
et tel que cette propriété soit vérifiée pour tout entier $\leq s'$. Comme $V=W_0$, la construction de $\tilde{Z}$ garantit que $s' \geq 1$ (la perte de composantes au premier cran ne peut être plus grande que la perte totale). 

Cette clause technique tient au fait que la variété obstructrice dans le lemme de zéros ne contient pas nécessairement toutes les composantes de $V + \Sigma_0$, mais une proportion importante par un principe de tiroirs. On a donc introduit le sous-ensemble $\Sigma$ de $\Sigma_0$ qu'on ne contrôle pas et (\ref{clause technique}) montre que la perte de composantes dans la descente est limitée. Elle sera compensée par les propriétés des indices d'obstruction en jeu.

\bigskip

Au cran suivant, on observe en particulier que: 
\begin{displaymath}
W_{s'} + \Sigma^{s'+1} \nsubseteq \tilde{Z},
\end{displaymath}
et on peut supposer sans perte de généralité que: 
\begin{displaymath}
W_{s'} \nsubseteq \tilde{Z}. 
\end{displaymath}
Par construction, il existe un sous-ensemble $\Gamma^{s'}$ de $\Sigma^{s'}$ vérifiant (\ref{clause technique}). Les deux fermés $W_{s'}$ et $\tilde{Z}$ contiennent alors $W_{s'-1} + \Gamma^{s'}$ (voir la définition \ref{7}).
Leur intersection $W_{s'} \cap \tilde{Z}$ est un fermé équidimensionnel de dimension $\mathrm{dim}(W_{s'}) -1$ et contient encore le fermé $W_{s'-1}+ \Gamma^{s'}$. De plus, son degré est majoré par le théorème de Bézout (quitte à sommer les degrés sur toutes les composantes irréductibles de $\tilde{Z}$): 
\begin{displaymath}
\de \big( W_{s'} \cap \tilde{Z} \big) \leq \de (W_{s'}) \de (\tilde{Z}).
\end{displaymath}
On choisit alors, comme dans le premier cas, une composante irréductible $Z_{s'}$ de cette intersection ``meilleure que la moyenne'', c'est-à-dire telle qu'il existe $\tilde{\Sigma}^{s'}$ avec: 
\begin{displaymath}
W_{s'-1} + \tilde{\Sigma}^{s'} \subset Z_{s'},
\end{displaymath}
et qu'on ait l'inégalité suivante: 
\begin{displaymath}
\frac{\de(Z_{s'})}{\de \big( W_{s'} \cap \tilde{Z} \big)} \leq \frac{\alpha_{W_{s'-1}}(\tilde{\Sigma}^{s'})}{\alpha_{W_{s'-1}}(\Gamma^{s'})}.
\end{displaymath}
On en déduit:
\begin{displaymath}
\frac{\alpha_{W_{s'-1}}(\Lambda^{s'})}{\alpha_{W_{s'-1}}(\tilde{\Sigma}^{s'})} \de(Z_{s'}) \leq  \frac{\alpha_{W_{s'-1}}(\Lambda^{s'})}{\alpha_{W_{s'-1}}(\Sigma^{s'})} \de(W_{s'}) \left( \Delta^3 p_{s+1} \frac{ \alpha_V(\Sigma_0)}{\alpha_V(\Sigma)} \omega(X)\right).
\end{displaymath}
A ce stade, on utilise les remarques précédant (\ref{recurrence obstruction}): 
\begin{displaymath}
\frac{\alpha_V(\Sigma_0)}{\alpha_V(\Sigma)} \omega(X) \leq \omega(V, \Sigma_0).
\end{displaymath}
De plus, soit $Z$ contenant un certain $V + \Lambda$ et réalisant $\omega(V, \Sigma^1 + \cdots + \Sigma^{s'})$. On a alors, quitte à considérer un translaté de $Z$ et en faisant une majoration grossière:
\begin{displaymath}
\omega(V, \Sigma_0) \leq \frac{\alpha_V(\Sigma_0)}{\alpha_V(\Lambda)} \de(Z) \leq p_{s'+1}^g \cdots p_s^g \omega(V, \Sigma^1 + \cdots + \Sigma^{s'}).
\end{displaymath}
Il reste enfin à utiliser la relation de récurrence de la définition \ref{def 7}:
\begin{displaymath}
\omega(V, \Sigma^1 + \cdots + \Sigma^{s'}) \leq \Delta^3 p_{s'} \omega(V_{s'-1}).
\end{displaymath}
En mettant bout à bout ces inégalités, il vient: 
\begin{displaymath}
\frac{\alpha_{W_{s'-1}}(\Lambda^{s'})}{\alpha_{W_{s'-1}}(\tilde{\Sigma}^{s'})} \de(Z_{s'}) \leq  \frac{\alpha_{W_{s'-1}}(\Lambda^{s'})}{\alpha_{W_{s'-1}}(\Sigma^{s'})} \de(W_{s'}) p_{s'} \omega(V_{s'-1}) \Big(  \Delta^6 p_{s'+1}^g \cdots p_{s+1}^g \Big) .
\end{displaymath}
\end{dem}

\bigskip

La proposition \ref{32} sera entièrement prouvée si on parvient à démontrer que la nouvelle suite $(s', \textbf{l}', {\bf \Sigma}', \textbf{W}')$ construite dans le lemme précédent est dans $\mathcal{W}(V)$:
\begin{lem}
\label{lemme final}
La suite $(s', \textbf{l}', {\bf \Sigma}', \textbf{W}')$ est dans $\mathcal{W}(V)$.
\end{lem}
\begin{dem}
Rappelons que la nouvelle suite est définie par: 
\begin{displaymath}
{\bf \Sigma}'= (\Sigma^1, \ldots, \Sigma^{s'-1}, \tilde{\Sigma}^{s'}),
\end{displaymath}
et: 
\begin{displaymath}
\textbf{W}'= (W_0, \ldots, W_{s'-1}, Z_{s'}).
\end{displaymath}
Le premier point est toujours vérifié par construction. 

Démontrons le second point pour $i=s'$, les autres cas étant inchangés. Si $s' \neq s+1$, on a:
\begin{eqnarray*}
\frac{\alpha_{W_{s'-1}}(\Lambda^{s'})}{\alpha_{W_{s'-1}}(\tilde{\Sigma}^{s'})} \de(Z_{s'}) & \leq & \frac{\alpha_{W_{s'-1}}(\Lambda^{s'})}{\alpha_{W_{s'-1}}(\Sigma^{s'})} \de(W_{s'}) p_{s'} \omega(V_{s'-1}) \Delta^6 (p_{s'+1} \cdots p_{s+1})^{g}  \\
& \leq &  \Big(  p_{s'} \omega(V_{s'-1}) \Delta^6 (P_{s'+1} \cdots P_k)^{g}  \Big)^{\mathrm{codim}(W_{s'})}  p_{s'} \omega(V_{s'-1}) \Delta^6 (p_{s'+1} \cdots p_{s+1})^{g}  \\
& \leq & \Big( p_{s'} \omega(l_{s'}, V_{s'-1}) \Delta^6  (P_{s'+1} \cdots P_{k})^{g} \Big)^{\mathrm{codim}(Z_{s'})},
\end{eqnarray*}
puisque dans ce cas: codim$(Z_{s'})=$ codim$(W_{s'})+1$. Si $s'=s+1$, la majoration du lemme \ref{33} donne directement le résultat.

Enfin, la ``relation de récurrence'' entre indices d'obstruction est trivialement réalisée pour $i \neq s'$; pour $i=s'=s+1$, c'est l'inégalité (\ref{recurrence obstruction}); et pour $i=s'$ et $s' \neq s+1$, il n'y a rien à montrer. 
\end{dem}

\bigskip

\noindent
\textbf{Preuve} (du théorème \ref{5}) 

On peut maintenant démontrer le théorème \ref{5}. Soit $V$ une sous-variété propre de $A$ qui n'est pas incluse dans un translaté de sous-variété abélienne de $A$. Par les propositions \ref{31} et \ref{32}, la variété $V$ contredit l'inégalité suivante: 
\begin{displaymath}
\me(V) \ \omega(V) < \frac{1}{\Delta^{\left(5g(k+1)\right)^{k+1}}},
\end{displaymath}
avec $\Delta=C_0^2 \mathrm{log}(3 \de(V))$. On en déduit: 
\begin{displaymath}
\me(V) \geq \frac{C(A)}{\omega(V)} \times
(\mathrm{log}(3 \de(V)))^{-\lambda(k)},
\end{displaymath}
où $\lambda(k)=(5g(k+1))^{k+1}$ et $C(A)= \frac{1}{C_0^{2 \lambda(g)}}$, qui ne dépend que de $A$.
{\begin{flushright}$\Box$\end{flushright}}

\bibliographystyle{fralpha}
\bibliography{biblio}

\begin{thebibliography}{MvdG07}
\expandafter\ifx\csname fonteauteurs\endcsname\relax
\def\fonteauteurs{\scshape}\fi

\bibitem[AD99]{Amoroso-David99}
F.~\bgroup\fonteauteurs\bgroup Amoroso\egroup\egroup{} et
  S.~\bgroup\fonteauteurs\bgroup David\egroup\egroup{} :
\newblock Le probl{\`e}me de {L}ehmer en dimension sup{\'e}rieure.
\newblock {\em J. reine angew. Math.}, 513\string:\penalty500\relax 145--179,
  1999.

\bibitem[AD03]{Amoroso-David03}
F.~\bgroup\fonteauteurs\bgroup Amoroso\egroup\egroup{} et
  S.~\bgroup\fonteauteurs\bgroup David\egroup\egroup{} :
\newblock Minoration de la hauteur normalis{\'e}e dans un tore.
\newblock {\em J. Inst. Math. Jussieu}, 2(3)\string:\penalty500\relax 335--381,
  2003.

\bibitem[AD07]{Amoroso-David07}
F.~\bgroup\fonteauteurs\bgroup Amoroso\egroup\egroup{} et
  S.~\bgroup\fonteauteurs\bgroup David\egroup\egroup{} :
\newblock Un th{\'e}or{\`e}me de z{\'e}ros dans les groupes alg{\'e}briques
  commutatifs.
\newblock {\em Manuscript en pr{\'e}paration}, 2007.

\bibitem[Amo07]{Amoroso07}
F~\bgroup\fonteauteurs\bgroup Amoroso\egroup\egroup{} :
\newblock Bogomolov on tori revisited.
\newblock {\em Pr{\'e}publication}, 2007.

\bibitem[BGS94]{Bost-Gillet-Soule94}
J.B. \bgroup\fonteauteurs\bgroup Bost\egroup\egroup{},
  H.~\bgroup\fonteauteurs\bgroup Gillet\egroup\egroup{} et
  C.~\bgroup\fonteauteurs\bgroup Soul{\'e}\egroup\egroup{} :
\newblock Heights of projective varieties.
\newblock {\em Journal of the A.M.S.}, 7(4)\string:\penalty500\relax 903--1027,
  1994.

\bibitem[BK06]{Bost-Kunnemann06}
J.B. \bgroup\fonteauteurs\bgroup Bost\egroup\egroup{} et
  K.~\bgroup\fonteauteurs\bgroup Künnemann\egroup\egroup{} :
\newblock Hermitian vector bundles and extension groups on arithmetic
  varieties.
\newblock {\em Pr{\'e}publication}, 2006.

\bibitem[BMZ99]{Bombieri-Masser-Zannier99}
E.~\bgroup\fonteauteurs\bgroup Bombieri\egroup\egroup{},
  D.~\bgroup\fonteauteurs\bgroup Masser\egroup\egroup{} et
  U.~\bgroup\fonteauteurs\bgroup Zannier\egroup\egroup{} :
\newblock Intersecting a curve with algebraic subgroups of multiplicative
  groups.
\newblock {\em Internat. Math. Res. Notices}, 20\string:\penalty500\relax
  1119--1140, 1999.

\bibitem[Bos96a]{Bost96}
J.B. \bgroup\fonteauteurs\bgroup Bost\egroup\egroup{} :
\newblock Intrinsic heights of stable varieties and abelian varieties.
\newblock {\em Duke Math. J.}, 82(1)\string:\penalty500\relax 21--70, 1996.

\bibitem[Bos96b]{Bost95}
J.B. \bgroup\fonteauteurs\bgroup Bost\egroup\egroup{} :
\newblock P{\'e}riodes et isog{\'e}nies des vari{\'e}t{\'e}s ab{\'e}liennes sur
  les corps de nombres (d'apr{\`e}s {D}. {M}asser et {G}. {W}üstholz).
\newblock {\em Séminaire Bourbaki, Astérisque}, 237\string:\penalty500\relax
  115--161, 1996.

\bibitem[Bos01]{Bost01}
J.B. \bgroup\fonteauteurs\bgroup Bost\egroup\egroup{} :
\newblock Algebraic leaves of algebraic foliations over number fields.
\newblock {\em Publ. Math. IH{\'E}S}, 93\string:\penalty500\relax 161--221,
  2001.

\bibitem[BZ95]{Bombieri-Zannier95}
E.~\bgroup\fonteauteurs\bgroup Bombieri\egroup\egroup{} et
  U.~\bgroup\fonteauteurs\bgroup Zannier\egroup\egroup{} :
\newblock Algebraic points on subvarieties of $\mathbb{G}_m^n$.
\newblock {\em Internat. Math. Res. Notices}, 7\string:\penalty500\relax
  333--347, 1995.

\bibitem[BZ96]{Bombieri-Zannier96}
E.~\bgroup\fonteauteurs\bgroup Bombieri\egroup\egroup{} et
  U.~\bgroup\fonteauteurs\bgroup Zannier\egroup\egroup{} :
\newblock Heights of algebraic points on subvarieties of abelian varieties.
\newblock {\em Ann. Scuola Norm. Pisa}, 23(4)\string:\penalty500\relax
  779--792, 1996.

\bibitem[Cha88]{Chardin89}
M.~\bgroup\fonteauteurs\bgroup Chardin\egroup\egroup{} :
\newblock Une majoration de la fonction de {H}ilbert et ses cons{\'e}quences
  pour l'interpolation alg{\'e}brique.
\newblock {\em Bull. Soc. math. France}, 117\string:\penalty500\relax 305--318,
  1988.

\bibitem[Cha90]{Chardin90}
M.~\bgroup\fonteauteurs\bgroup Chardin\egroup\egroup{} :
\newblock Contributions {\`a} l'alg{\`e}bre commutative effective et {\`a} la
  th{\'e}orie de l'{\'e}limination.
\newblock {\em Th{\`e}se de {D}octorat, Universit{\'e} Paris {V}{I}}, 1990.

\bibitem[Che06]{Chen06}
H.~\bgroup\fonteauteurs\bgroup Chen\egroup\egroup{} :
\newblock Positivit{\'e} en g{\'e}om{\'e}trie alg{\'e}brique et en
  g{\'e}om{\'e}trie d'{A}rakelov.
\newblock {\em Th{\`e}se de {D}octorat, {U}niversit{\'e} {P}aris {X}{I}}, 2006.

\bibitem[Dav91]{David91}
S.~\bgroup\fonteauteurs\bgroup David\egroup\egroup{} :
\newblock Fonctions th{\^e}ta et points de torsion des vari{\'e}t{\'e}s
  ab{\'e}liennes.
\newblock {\em Compositio Math.}, 78(2)\string:\penalty500\relax 121--160,
  1991.

\bibitem[DH00]{David-Hindry00}
S.~\bgroup\fonteauteurs\bgroup David\egroup\egroup{} et
  M.~\bgroup\fonteauteurs\bgroup Hindry\egroup\egroup{} :
\newblock Minoration de la hauteur de {N}{\'e}ron-{T}ate sur les
  vari{\'e}t{\'e}s ab{\'e}liennes de type {C}.{M}.
\newblock {\em J. reine angew. Math.}, 529\string:\penalty500\relax 1--74,
  2000.

\bibitem[DP00]{David-Philippon00}
S.~\bgroup\fonteauteurs\bgroup David\egroup\egroup{} et
  P.~\bgroup\fonteauteurs\bgroup Philippon\egroup\egroup{} :
\newblock Sous-vari{\'e}t{\'e}s de torsion des vari{\'e}t{\'e}s
  semi-ab{\'e}liennes.
\newblock {\em C. R. Acad. Sci. Paris}, 331(1)\string:\penalty500\relax
  587--592, 2000.

\bibitem[DP02]{David-Philippon02}
S.~\bgroup\fonteauteurs\bgroup David\egroup\egroup{} et
  P.~\bgroup\fonteauteurs\bgroup Philippon\egroup\egroup{} :
\newblock Minoration des hauteurs normalis{\'e}es des sous-vari{\'e}t{\'e}s de
  vari{\'e}t{\'e}s ab{\'e}liennes {I}{I}.
\newblock {\em Comment. Math. Helv.}, 77\string:\penalty500\relax 639--700,
  2002.

\bibitem[Ful84]{Fulton84}
W.~\bgroup\fonteauteurs\bgroup Fulton\egroup\egroup{} :
\newblock {\em Intersection theory}.
\newblock Springer-Verlag, Berlin, 1984.

\bibitem[Gal07]{Galateau07}
A.~\bgroup\fonteauteurs\bgroup Galateau\egroup\egroup{} :
\newblock Probl{\`e}me de {B}ogomolov sur les vari{\'e}t{\'e}s ab{\'e}liennes.
\newblock {\em Th{\`e}se de {D}octorat, {U}niversit{\'e} {P}aris {V}{I}}, 2007.

\bibitem[Gau06]{Gaudron06}
{\'E}.~\bgroup\fonteauteurs\bgroup Gaudron\egroup\egroup{} :
\newblock Formes lin{\'e}aires de logarithmes effectives sur les
  vari{\'e}t{\'e}s ab{\'e}liennes.
\newblock {\em Annales de l'{\'E}NS}, 2006.

\bibitem[Gau07]{Gaudron07}
{\'E}.~\bgroup\fonteauteurs\bgroup Gaudron\egroup\egroup{} :
\newblock Pentes des fibr{\'e}s vectoriels ad{\'e}liques sur un corps global.
\newblock {\em Rendiconti di Padova}, 2007.

\bibitem[GH78]{Griffiths-Harris78}
P.~\bgroup\fonteauteurs\bgroup Griffiths\egroup\egroup{} et
  J.~\bgroup\fonteauteurs\bgroup Harris\egroup\egroup{} :
\newblock {\em Principles of algebraic geometry}.
\newblock John Wiley and sons, New-York, 1978.

\bibitem[Gra01]{Graftieaux01}
P.~\bgroup\fonteauteurs\bgroup Graftieaux\egroup\egroup{} :
\newblock Formal groups and {I}sogeny theorem.
\newblock {\em Duke Math. J.}, 106\string:\penalty500\relax 81--121, 2001.

\bibitem[GS92]{Gillet-Soule92}
H.~\bgroup\fonteauteurs\bgroup Gillet\egroup\egroup{} et
  C.~\bgroup\fonteauteurs\bgroup Soul{\'e}\egroup\egroup{} :
\newblock An arithmetic {R}iemann-{R}och theorem.
\newblock {\em Invent. Math.}, 110(3)\string:\penalty500\relax 473--543, 1992.

\bibitem[Hab06]{Habegger06}
P.~\bgroup\fonteauteurs\bgroup Habegger\egroup\egroup{} :
\newblock A {B}ogomolov property modulo algebraic subgroups.
\newblock {\em Pr{\'e}publication}, 2006.

\bibitem[Har77]{Hartshorne77}
R.~\bgroup\fonteauteurs\bgroup Hartshorne\egroup\egroup{} :
\newblock {\em Algebraic {G}eometry}, volume~52 de {\em Graduate Texts in
  Mathematics}.
\newblock Springer-Verlag, New York, 1977.

\bibitem[Hin88]{Hindry88}
M.~\bgroup\fonteauteurs\bgroup Hindry\egroup\egroup{} :
\newblock Autour d'une conjecture de {S}. {L}ang.
\newblock {\em Invent. Math.}, 94\string:\penalty500\relax 575--603, 1988.

\bibitem[HS00]{Hindry-Silverman00}
M.~\bgroup\fonteauteurs\bgroup Hindry\egroup\egroup{} et
  J.~\bgroup\fonteauteurs\bgroup Silverman\egroup\egroup{} :
\newblock {\em {D}iophantine {G}eometry: {A}n {I}ntroduction}, volume 201 de
  {\em Graduate Texts in Mathematics}.
\newblock Springer-Verlag, New York, 2000.

\bibitem[IR80]{Ireland-Rosen80}
K.~\bgroup\fonteauteurs\bgroup Ireland\egroup\egroup{} et
  M.~\bgroup\fonteauteurs\bgroup Rosen\egroup\egroup{} :
\newblock {\em A {C}lassical {I}ntroduction to {M}odern {N}umber {T}heory},
  volume~84 de {\em Graduate Texts in Mathematics}.
\newblock Springer-Verlag, New York, 1980.

\bibitem[Kob75]{Koblitz75}
N.~\bgroup\fonteauteurs\bgroup Koblitz\egroup\egroup{} :
\newblock $p$-adic variation of the zeta-function over families of varieties
  defined over finite fields.
\newblock {\em Compositio Math.}, 31(2), 1975.

\bibitem[Lan86]{Lang86}
S.~\bgroup\fonteauteurs\bgroup Lang\egroup\egroup{} :
\newblock {\em Algebraic number theory}.
\newblock Springer-Verlag, 1986.

\bibitem[Lan02]{Lang02}
S.~\bgroup\fonteauteurs\bgroup Lang\egroup\egroup{} :
\newblock {\em Algebra}, volume 211 de {\em Graduate Texts in Mathematics}.
\newblock Springer-Verlag, New York, 2002.

\bibitem[Lau83]{Laurent83}
M.~\bgroup\fonteauteurs\bgroup Laurent\egroup\egroup{} :
\newblock Minoration de la hauteur de {N}{\'e}ron-{T}ate.
\newblock {\em S{\'e}minaire de th{\'e}orie des nombres de Paris, 1981-1982,
  Progr. Math.}, 38\string:\penalty500\relax 137--152, 1983.

\bibitem[LR85]{Lange-Ruppert85}
H.~\bgroup\fonteauteurs\bgroup Lange\egroup\egroup{} et
  W.~\bgroup\fonteauteurs\bgroup Ruppert\egroup\egroup{} :
\newblock Complete systems of addition laws on abelian varieties.
\newblock {\em Invent. Math.}, 79(3)\string:\penalty500\relax 603--610, 1985.

\bibitem[Mau07]{Maurin07}
G.~\bgroup\fonteauteurs\bgroup Maurin\egroup\egroup{} :
\newblock Conjecture de {Z}ilber-{P}ink pour les courbes trac{\'e}es sur des
  tores.
\newblock {\em Pr{\'e}publication de l'{I}nstitut {F}ourier}, 696, 2007.

\bibitem[MB90]{Moret-Bailly90}
L.~\bgroup\fonteauteurs\bgroup Moret-Bailly\egroup\egroup{} :
\newblock Sur l'{\'e}quation fonctionnelle de la fonction thêta de {R}iemann.
\newblock {\em Compositio Math.}, 75\string:\penalty500\relax 203--217, 1990.

\bibitem[Mou04]{Mourougane04}
C.~\bgroup\fonteauteurs\bgroup Mourougane\egroup\egroup{} :
\newblock Computations of {B}ott-{C}hern classes on $\mathbb{P}({E})$.
\newblock {\em Duke Math. J.}, 24(2)\string:\penalty500\relax 389--420, 2004.

\bibitem[Mum74]{Mumford}
D.~\bgroup\fonteauteurs\bgroup Mumford\egroup\egroup{} :
\newblock {\em Abelian {V}arieties}.
\newblock Tata Lecture Notes. Cambridge University Press, 1974.

\bibitem[MvdG07]{Moonen}
B.~\bgroup\fonteauteurs\bgroup Moonen\egroup\egroup{} et G.~van~der
  \bgroup\fonteauteurs\bgroup Geer\egroup\egroup{} :
\newblock {\em Abelian varieties}.
\newblock Version pr{\'e}liminaire.
  http://staff.science.uva.nl/~bmoonen/boek/BookAV.html, 2007.

\bibitem[Noo95]{Noot95}
R.~\bgroup\fonteauteurs\bgroup Noot\egroup\egroup{} :
\newblock Abelian varieties - {G}alois representations and properties of
  ordinary reduction.
\newblock {\em Compositio Math.}, 97\string:\penalty500\relax 161--171, 1995.

\bibitem[Ogu82]{Ogus83}
A.~\bgroup\fonteauteurs\bgroup Ogus\egroup\egroup{} :
\newblock Hodge cycles and crystalline cohomology.
\newblock {\em Lecture Notes in Mathematics}, 900, 1982.

\bibitem[Phi95]{Philippon95}
P.~\bgroup\fonteauteurs\bgroup Philippon\egroup\egroup{} :
\newblock Sur des hauteurs alternatives {I}{I}{I}.
\newblock {\em J. Math. Pures Appl.}, 74(4)\string:\penalty500\relax 345--365,
  1995.

\bibitem[Pin98]{Pink98}
R.~\bgroup\fonteauteurs\bgroup Pink\egroup\egroup{} :
\newblock $l$-adic algebraic monodromy groups, cocharacters, and the
  {M}umford-{T}ate conjecture.
\newblock {\em J. reine angew. Math.}, 495\string:\penalty500\relax 187--237,
  1998.

\bibitem[Pin04]{Pink-online}
R.~\bgroup\fonteauteurs\bgroup Pink\egroup\egroup{} :
\newblock {\em Finite {G}roup {S}chemes}.
\newblock {N}otes de {C}ours \\
  http://www.math.ethz.ch/~pink/FiniteGroupSchemes.html, 2004.

\bibitem[Pin05]{Pink05}
R.~\bgroup\fonteauteurs\bgroup Pink\egroup\egroup{} :
\newblock A common generalization of the conjectures of {A}ndr{\'e}-{O}ort,
  {M}anin-{M}umford and {M}ordell-{L}ang.
\newblock {\em Pr{\'e}publication}, 2005.

\bibitem[Ray74]{Raynaud74}
M.~\bgroup\fonteauteurs\bgroup Raynaud\egroup\egroup{} :
\newblock Sch{\'e}mas en groupes de type $(p, \ldots, p)$.
\newblock {\em Bulletin de la S.M.F.}, 2\string:\penalty500\relax 241--280,
  1974.

\bibitem[Ray83]{Raynaud83}
M.~\bgroup\fonteauteurs\bgroup Raynaud\egroup\egroup{} :
\newblock Courbes sur une vari{\'e}t{\'e} ab{\'e}lienne et points de torsion.
\newblock {\em Invent. Math.}, 71\string:\penalty500\relax 207--233, 1983.

\bibitem[Sam03]{Samuel}
P.~\bgroup\fonteauteurs\bgroup Samuel\egroup\egroup{} :
\newblock {\em Th{\'e}orie alg{\'e}brique des nombres}.
\newblock Hermann, Paris, 2003.

\bibitem[Ser68]{Serre68}
J.P. \bgroup\fonteauteurs\bgroup Serre\egroup\egroup{} :
\newblock {\em Abelian l-adic representations and elliptic curves}.
\newblock Benjamin, New York, 1968.

\bibitem[Ten95]{Tenenbaum}
G.~\bgroup\fonteauteurs\bgroup Tenenbaum\egroup\egroup{} :
\newblock {\em Introduction {\`a} la th{\'e}orie analytique et probabiliste des
  nombres}.
\newblock Collection S.M.F., Paris, 1995.

\bibitem[Ull98]{Ullmo96}
E.~\bgroup\fonteauteurs\bgroup Ullmo\egroup\egroup{} :
\newblock Positivit{\'e} et discr{\'e}tion des points alg{\'e}briques des
  courbes.
\newblock {\em Ann. of Math.}, 147\string:\penalty500\relax 167--179, 1998.

\bibitem[Zha92]{Zhang92}
S.~\bgroup\fonteauteurs\bgroup Zhang\egroup\egroup{} :
\newblock Positive line bundles on arithmetic surfaces.
\newblock {\em Ann. of Math.}, 136\string:\penalty500\relax 569--587, 1992.

\bibitem[Zha95a]{Zhang-Positive95}
S.~\bgroup\fonteauteurs\bgroup Zhang\egroup\egroup{} :
\newblock Positive line bundles on arithmetic varieties.
\newblock {\em J. Amer. Math. Soc.}, 8\string:\penalty500\relax 187--221, 1995.

\bibitem[Zha95b]{Zhang95}
S.~\bgroup\fonteauteurs\bgroup Zhang\egroup\egroup{} :
\newblock Small points and adelic metrics.
\newblock {\em J. Algebraic Geom.}, 4\string:\penalty500\relax 281--300, 1995.

\bibitem[Zha98]{Zhang98}
S.~\bgroup\fonteauteurs\bgroup Zhang\egroup\egroup{} :
\newblock Equidistribution of small points on abelian varieties.
\newblock {\em Ann. of Math.}, 147\string:\penalty500\relax 159--165, 1998.

\bibitem[Zil02]{Zilber02}
B.~\bgroup\fonteauteurs\bgroup Zilber\egroup\egroup{} :
\newblock Intersecting varieties with tori.
\newblock {\em J. London Math. Soc.}, 65\string:\penalty500\relax 27--44, 2002.

\end{thebibliography}

\end{document}